\documentclass[a4paper]{article}
\usepackage[affil-it]{authblk}

\usepackage{array}
\usepackage{color}
\usepackage{tabularx}
\usepackage{graphicx}
\usepackage{amsfonts}	
\usepackage{moreverb}
\usepackage{dsfont}
\usepackage{grffile}
\usepackage{bm}
\usepackage{multirow}
\usepackage{soul}
\usepackage{booktabs}
\newcommand\BibTeX{{\rmfamily B\kern-.05em \textsc{i\kern-.025em b}\kern-.08em
T\kern-.1667em\lower.7ex\hbox{E}\kern-.125emX}}

\newcommand{\Int}{\displaystyle\int}
\newcommand{\Sum}{\displaystyle\sum}

\renewcommand{\v}{\mathbf{v}}

\newcommand{\be}{\begin{equation}}
\newcommand{\ee}{\end{equation}}
\newcommand{\bdm}{\begin{displaymath}}
\newcommand{\edm}{\end{displaymath}}
\renewcommand{\v}{\mathbf{v}}

\newcommand{\apriori}{\textit{a priori} }
\newcommand{\aposteriori}{\textit{a posteriori} }

\newfont{\numerikEleven}{ecrm1000}
\newfont{\numerikTen}{cmss10}
\newfont{\numerikNine}{cmss9}
\newfont{\numerikEight}{cmss8}
\newfont{\numerikSeven}{cmss7}
\newfont{\numerikSix}{cmss6}

\usepackage{amsmath,amsfonts,float,amssymb,cases}
\usepackage[amsmath,thmmarks]{ntheorem}

\theoremseparator{}
\theorembodyfont{\normalfont}
\theoremsymbol{}
\newtheorem{teoI}{Theorem}%[section]

\newtheorem{remI}[teoI]{Remark}

\theoremstyle{nonumberplain}
\theoremheaderfont{\normalfont\bfseries}
\theoremseparator{.}
\theoremsymbol{\ensuremath{_\blacksquare}}

\definecolor{violet}{rgb}{0.93, 0.51, 0.93}
\definecolor{indigo}{rgb}{0.29, 0, 0.51}
\definecolor{blue}{rgb}{0, 0, 1}
\definecolor{green}{rgb}{0, 0.50, 0}
\definecolor{yellow}{rgb}{1, 1, 0}
\definecolor{orange}{rgb}{1, 0.65, 0}
\definecolor{red}{rgb}{1, 0, 0}

\usepackage{tikz}
\usepackage[export]{adjustbox}
\usepackage{pgfplots}
\pgfplotsset{compat = newest}
\newcommand{\exactsol}{\raisebox{0pt}{\tikz{\draw[red,line width = 1.0pt](0.,0.8mm) -- (5.5mm,0.8mm);}}}
\newcommand{\rhosup}{\raisebox{0pt}{\tikz{\draw[green,line width = 1.0pt](0.,0.8mm) -- (5.5mm,0.8mm);}}}
\newcommand{\rhosub}{\raisebox{0pt}{\tikz{\draw[blue,line width = 1.0pt](0.,0.8mm) -- (5.5mm,0.8mm);}}}

\newcommand{\blackblack}{\tikz{\draw[fill,black] (2.75mm,0.75mm) circle (0.65mm);}}
\newcommand{\redred}{\tikz{\draw[fill,red] (2.75mm,0.75mm) circle (0.65mm);}}
\newcommand{\greengreen}{\tikz{\draw[fill,green] (2.75mm,0.75mm) circle (0.65mm);}}

\newcommand{\CPDviolet}{\tikz{\draw[fill=violet] (2.75mm,0.75mm) circle (0.95mm);}}
\newcommand{\CPDindigo}{\tikz{\draw[fill={rgb:red,0.29; green,0; blue,0.51}] (2.75mm,0.75mm) circle (0.95mm);}}
\newcommand{\CPDblue}{\tikz{\draw[fill=blue] (2.75mm,0.75mm) circle (0.95mm);}}
\newcommand{\CPDgreen}{\tikz{\draw[fill=green] (2.75mm,0.75mm) circle (0.95mm);}}
\newcommand{\CPDyellow}{\tikz{\draw[fill=yellow] (2.75mm,0.75mm) circle (0.95mm);}}
\newcommand{\CPDorange}{\tikz{\draw[fill=orange] (2.75mm,0.75mm) circle (0.95mm);}}
\newcommand{\CPDred}{\tikz{\draw[fill=red] (2.75mm,0.75mm) circle (0.95mm);}}
\newcommand{\CPDwhite}{\tikz{\draw[fill=white] (2.75mm,0.75mm) circle (0.95mm);}}

\newcommand{\velocity}{\raisebox{2pt}{\tikz{\draw[black,line width = 1.0pt](0.,0.8mm) -- (5.5mm,0.8mm);}}}
\newcommand{\phiphi}{\raisebox{2pt}{\tikz{\draw[red,line width = 1.0pt](0.,0.8mm) -- (5.5mm,0.8mm);}}}

\newfont{\tablepolice}{cmss8}
\usepackage[linesnumbered,ruled,vlined]{algorithm2e}
\SetKwInput{KwData}{Input}
\SetKwInput{KwResult}{Output}

\usepackage{graphicx}
\usepackage{graphbox}

\begin{document} 
\title{\aposteriori stabilized sixth-order finite volume scheme\\with adaptive stencil construction\\ --- Basics for the 1D steady-state hyperbolic equations} 

\author{Gaspar~J.~Machado%
  \thanks{Electronic address: \texttt{gjm@math.uminho.pt}}}
\affil{Centre of Physics and Department of Mathematics, University of Minho\\
  Campus of Azur\'em, 4800-058 Guimar\~aes, Portugal}

\author{St{\'e}phane Clain%
  \thanks{Electronic address: \texttt{clain@math.uminho.pt}}}
\affil{Centre of Physics and Department of Mathematics, University of Minho\\
  Campus of Azur\'em, 4800-058 Guimar\~aes, Portugal}

\author{Rapha{\"e}l Loub{\`e}re%
  \thanks{Electronic address: \texttt{raphael.loubere@u-bordeaux.fr}}}
\affil{CNRS and Institut de Math\'{e}matiques de Bordeaux (IMB)\\
  Universit{\'e} de Bordeaux, Talence, France}

\date{September, 2020}

\maketitle

%-------------------------------------------------------
% ABSTRACT
\begin{abstract}
% \PACS{PACS code1 \and PACS code2 \and more}
% \subclass{MSC code1 \and MSC code2 \and more}
% 65M08, 65M12, 65M99, 65H10, 76-10
We propose an adaptive stencil construction for high order accurate finite volume schemes \aposteriori stabilized devoted to solve one-dimensional steady-state hyperbolic equations. 
High-accuracy (up to the sixth-order presently) is achieved thanks to polynomial reconstructions while stability is provided with an \aposteriori MOOD method which controls the cell polynomial degree for eliminating non-physical oscillations in the vicinity of discontinuities. 
We supplemented this scheme with a stencil construction 
%which avoids the presence of troubled cells, 
allowing to reduce even further the numerical dissipation.
The stencil is shifted away from troubles (shocks, discontinuities, etc.) leading to less oscillating polynomial reconstructions. % and more accuracy.
Experimented on linear, B\"urgers', and Euler equations, we demonstrate that the adaptive stencil technique manages to retrieve smooth solutions with optimal order of accuracy but also irregular ones without spurious oscillations. 
Moreover we numerically show that the approach allows to %drastically 
reduce the dissipation still maintaining the essentially non-oscillatory behavior.
%\keywords{Finite Volume \and MOOD \and Adaptive stencil \and Steady state solution \and Euler equations \and High-order}
\end{abstract}

\tableofcontents

%==================================================================================
%
% I N T R O D U C T I O N
%
%\section{Introduction} \label{sec:introduction}
\section{Introduction} \label{sec:introduction}

Numerical approximations of the steady solution of the Euler equations is a long lasting activity in computational physics since the 70's, see \cite{Pu}. It has led to the development of a tremendous amount of techniques over the years, for instance the implicit time algorithm \cite{BW1} with high-order finite difference methods \cite{BW2} or some of the first simulations for two- and three-dimensional complex geometries calculated on Illiac IV Computer \cite{St,PuSt,PuLo}. The steady-state solution was approached as the limit stage of a non-stationary problem using a fictitious time step, leading to the so-called ``time marching method'' \cite{VeMa}. 

Regardless of the numerical method in hand some antagonist requirements have to be fulfilled. Firstly, a high accuracy is expected in the regular zones of the solution demanding a high-accurate numerical method without (or with an extremely) low numerical dissipation.
Secondly, a robust essentially non-oscillatory (ENO) solution, \textit{i.e.}, exempt from spurious oscillations, is needed where high gradients are encountered. These oscillations being a resultant from the Gibbs phenomenon they are generally damped by some numerical dissipation embedded in the numerical method. High-order methods supported by stabilization techniques to prevent non physical oscillations have been developed for more than half a century. In the Finite Volume (FV) context, second order methods date back to the 70's and, nowadays, very high order finite volume methods are well-established. While the polynomial reconstruction is a common ingredient to increase the accuracy, the limiting techniques still represent a challenging goal: to preserve the stability without sacrificing the accuracy. 

The crucial point lies in the detection of a potential discontinuity or steep gradient. More precisely, if one has to perform a polynomial reconstruction for a given cell $i$, a stencil $\mathcal{S}_i$ composed of neighbor cells is required and its associated values are interpolated by a high order polynomial. The resulting interpolation may have strong oscillations if some cells present large variations, that is larger than $\mathcal O(1/h)$ where $h$ is the characteristic mesh size parameter. Such unacceptable oscillations usually lead to a polynomial degree reduction in order to add numerical dissipation. To avoid such a situation, one has to carefully pick up the cells in the vicinity of cell $i$ to provide a stencil with only clean cells if possible, \text{i.e.}, cells away from a discontinuity. The principle of carefully selecting the stencil at least dates back to the very first versions of (Weighted)-Essentially-Non-Oscillatory (W)ENO  \cite{harten,liu,shu_efficient_weno}, where ``The basic idea is to avoid including the discontinuous cell in the stencil, if possible.'', see for instance page~342 of C.W.~Shu's monograph \cite{Shu1998}. In one respect, any ENO/WENO  algorithm produces some ``fuzzy'' stencil since it uses a non-linear combination between several polynomial reconstructions, each involving a different stencil. The resulting polynomial is not strictly speaking deriving from a unique stencil where the interpolation would be achieved, but a weighted combination of stencils \textit{a priori} deriving from smoothness indicators. 

In a recent work \cite{1D_SS} we have designed a new family of high order accurate Finite Volume (FV) schemes. High-accuracy is achieved thanks to polynomial reconstructions on centered stencils while stability and robustness are gained by an \aposteriori Multidimensional Optimal Order Detection (MOOD) method which controls the cell polynomial degree, eliminating non-physical oscillations in the vicinity of discontinuities by a reduction up to degree zero when needed. Such a procedure demands (i) to discriminate between troubled and valid cells, (ii) a cascade of decreasing polynomial degrees to be successively tested when oscillations are detected, and (iii) a parachute scheme corresponding to the last, viscous, and robust scheme. 

In \cite{1D_SS} numerical experimentation has been carried on advection, B\"urgers’, and Euler equations, and we have demonstrated that the scheme manages to retrieve smooth solutions with optimal order of accuracy but also irregular solutions without spurious oscillations. The \aposteriori MOOD loop allows to produce a valid solution by polynomial degree decrementing. Recall that a lower degree demands a narrower centered stencil for the reconstruction, which at last reduces to an empty one in the case of degree zero. As such this shrinking of the stencil width reduces \textit{de facto} the number of neighbors in the stencil, lowering the chances to ``see'' a discontinuity which may pollute the polynomial reconstruction. 

In the present paper, we propose an extension of the \aposteriori MOOD method assuming that the stencil is not given \textit{a priori} but derives from an analysis of the local numerical solution. We label the technique as ``Adaptive Stencil (AS) method'' since the stencil shape is also an unknown of the problem. In conjunction with the decrementing of the polynomial degree, a shifting of the reconstruction stencil away from any discontinuity is achieved. Doing so a cell in the vicinity of a steep gradient could be updated with a maximal accuracy (maximal polynomial degree) at the price of a careful choice of this non-centered stencil of neighbors used during the polynomial reconstruction.
To this end, we shall extend the \textit{a posteriori} paradigm to determine, cell by cell, the optimal stencil that provides both the accuracy and the stability. Roughly speaking, we produce a first candidate solution following the MOOD algorithm and by using centered stencils. We then obtain a map of the best polynomial degrees, called the Cell Polynomial Degree (CPD) map, that represents an approximation of the solution regularity. Based on this CPD map, we generate a new stencil for each cell by replacing the problematic cells, \textit{i.e.}, the ones having a low polynomial degree in the CPD map, by new neighborhood cells having higher polynomial degrees. In that way, we adapt the stencil in order to fit with the regular local areas
and try to sidestep the discontinuities. 

The goal of this work is to design such an AS algorithm to supplement the \aposteriori high-order MOOD FV scheme from \cite{1D_SS}. This proof of concept will be developed for the same three toy models: advection, B\"urgers’, and Euler equations in 1D for steady-state solutions for the sake of simplicity. Once designed, our AS algorithm will be tested against exact solutions to measure the gain in accuracy brought by such an approach. 
A Linear/Non-linear and Time marching iterative solvers will be considered in this work.

This paper is organized as follows. The governing equation and numerical schemes are first presented in section~\ref{sec:framework} and the classical MOOD scheme is recalled in section~\ref{sec:MOOD}. Then the core of this paper starts with the design of the AS method in section~\ref{sec:as} where the  algorithms are entirely described.
Then the numerical section~\ref{sec:numerics} gathers all evidences that the new approach outperforms the classical one on three 1D toy problems. At last conclusions and perspectives are drawn.

% end of file

%
% END  I N T R O D U C T I O N
%
%==================================================================================

%==================================================================================
%
% GOVERNING EQUATIONS AND SCHEMES
%
%\section{Governing equations and numerical schemes} \label{sec:equations_schemes}
\section{Governing equations and numerical schemes} \label{sec:framework} 

Let $\Omega=]x_{\text{L}},x_{\text{R}}[$ be the computational domain. We consider the steady state one-dimensional equation with source term
\begin{equation}
\label{eq:generic_equation}
\frac{\textrm{d}\mathbf{F}(\phi)}{\textrm{d}x}=\mathbf{S}, \quad\text{ in } \Omega,
\end{equation}
where $\mathbf{F}(\phi)$ stands for the physical flux and $\mathbf{S}\equiv\mathbf{S}(x)$ represents a regular source term. Dirichlet boundary conditions $\phi_{\text{L}}$ and $\phi_{\text{R}}$ are prescribed on the left and right boundaries of $\Omega$, respectively, which are only relevant at the inflow boundary interfaces. 

The system is assumed to be strictly hyperbolic and $\phi$ is the (scalar or vector) solution. For example, $\phi$ stands for a passive function for the advection problem while $\phi=u$ represents the scalar velocity for the B\"urgers equation. At last $\phi=(\rho,\rho u, E)$ represents the conservative variables of the Euler system of Partial Differential Equations (PDEs). In the same way, the source term is a scalar function or a vector valued function depending on the system in consideration.

\subsection{Mesh and notation} \label{ssec:formulation}

\begin{figure}%[ht]
\centering
\includegraphics{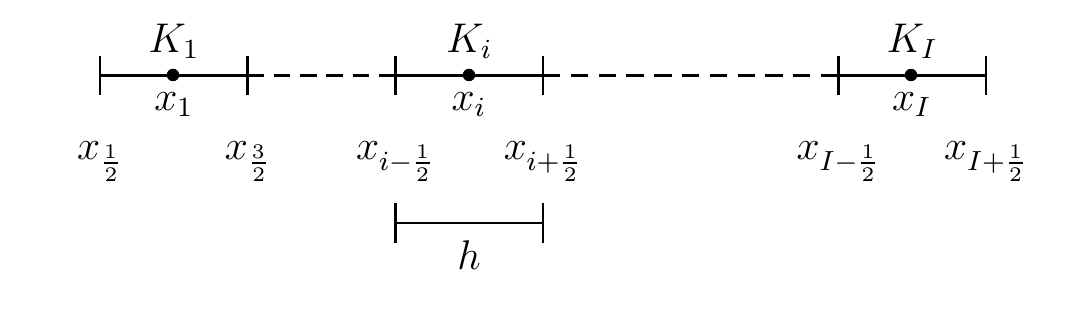}
\caption{Notation: points/interfaces and cells of the mesh.}
\label{figMeshNotation}
\end{figure}

The computational domain is discretized into a regular mesh $\mathcal T_h$ constituted of cells 
$K_i=[x_{i-\frac{1}{2}},x_{i+\frac{1}{2}}]$ for $i=1,\ldots,I$ (see Fig.~\ref{figMeshNotation}) of constant length 
$h=(x_\text{R}-x_\text{L})/I$ with center $x_i=\frac12(x_{i+\frac{1}{2}}+x_{i-\frac{1}{2}})$.
The boundaries of $\Omega$ are denoted by $x_{\frac{1}{2}}=x_{\text{L}}$ and $x_{I+\frac{1}{2}}=x_{\text{R}}$.
% ---- FIG -------
% ---- FIG -------
For any bounded function $\phi\equiv\phi(x)$ defined on $\Omega$, $\phi_i$ stands for an approximation of the mean value 
$\phi_i^{\text{ex}}$ of function $\phi$ over cell $K_i$, \textit{i.e.},
\[
\phi_i \approx \phi_i^{\text{ex}}=\frac{1}{h} \Int_{x_{i-\frac{1}{2}}}^{x_{i+\frac{1}{2}}} \phi(x) \, \textrm{d}x,
\qquad
i=1,\ldots,I,
\]
while vector $\Phi=\left(\phi_i\right)_{i=1,\ldots,I}\in\mathbb R^I$ gathers the values of these approximations.

\subsection{Polynomial reconstruction}\label{sec::def_stencil}

Let $K_i$ be a generic cell and $d\in\mathbb N$. A stencil $S_i$ is composed of $n_i$ neighbor cells (excluding cell $K_i$) such that $n_i>d$. The polynomial ${\phi}_{i}(x;d,S_i)$ of degree $d$ associated to cell $K_i$ is defined as
\[
{\phi}_{i}(x;d,S_i)=\phi_i+\sum_{k=1}^d R_{i,k}\left((x-x_i)^k -X_{i,k}\right),
\]
with $\bm{R}_i=(R_{i,k})_{k=1,\ldots,d}$ the vector gathering the unknown polynomial coefficients.
We set $X_{i,k}=\frac{1}{h}\int_{K_i} (x-x_i)^k \,\textrm{d}x$ to achieve the conservative property
\begin{equation}
\label{eq:conservation}
\frac{1}{h}\Int_{K_i} {\phi}_{i}(x;d,S_i)\,\textrm{d}x=\phi_i.
\end{equation}
For a given stencil $S_i$, we consider the quadratic functional 
\[
E_i(\bm{R}_i)=\sum_{j\in S_i} \left(\frac{1}{h} \int_{K_j} {\phi}_{i}(x;d,S_i)\,\textrm{d}x -\phi_j \right)^2.
\]
We denote by $\widehat{\bm{R}}_i=(\widehat{R}_{i,k})_{k=1,\ldots,d}$ the unique vector which minimizes the quadratic
functional and set $\widehat{\phi}_{i}(x;d,S_i)$ the associated polynomial that corresponds to the best approximation in the least squares sense of the data in stencil $S_i$, \textit{i.e.},
\[
\widehat{\phi}_{i}(x;d,S_i)=\phi_i+\sum_{k=1}^d \widehat{R}_{i,k}\left((x-x_i)^k -X_{i,k}\right).
\]
Such a reconstruction satisfies the conservation property~\eqref{eq:conservation} by construction.

%\textcolor{red}{RAPH: I do not think this remark is useful, we have defined $\phi_i$ and $\widehat\phi$ so what is the point of insisting on this remark? I am in favor of removing it.}
%\begin{remI} We have adopted the following notations: $\phi_i$ represents an approximation of $\phi$ over cell $K_i$ while ${\widehat\phi}_{i}(x;d;S_i)$ is a generic polynomial function of degree $d$ associated to $K_i$ and stencil $S_i$ which satisfies the conservation property (\ref{eq:conservation}).
%\end{remI}

\subsection{Generic finite volume scheme}

Integration equation~\eqref{eq:generic_equation} over a cell $K_i$ yields
\[
\Int_{x_{i-\frac{1}{2}}}^{x_{i+\frac{1}{2}}} 
\frac{\textrm{d}\mathbf{F}(\phi(x))}{\textrm{d}x}  \; \textrm{d}x  = 
\Int_{x_{i-\frac{1}{2}}}^{x_{i+\frac{1}{2}}} 
\mathbf{S}(x) \; \textrm{d}x , 
\]
and integration by parts provides the flux differential expression over cell $K_i$
\[
\left(\mathbf{F}(\phi(x_{i+\frac{1}{2}})-\mathbf{F}(\phi(x_{i-\frac{1}{2}})\right)-
h \left(\frac{1}{h}\Int_{x_{i-\frac{1}{2}}}^{x_{i+\frac{1}{2}}} \mathbf{S}(x) \; \textrm{d}x\right)=0.
\]
We substitute the exact relation with the numerical scheme 
\[
%\begin{equation}
%\label{eq:generic_FV}
\left(\bm{\mathcal{F}}(\phi_{i+\frac{1}{2},-},\phi_{i+\frac{1}{2},+};x_{i+\frac{1}{2}})-
\bm{\mathcal{F}}(\phi_{i-\frac{1}{2},-},\phi_{i-\frac{1}{2},+};x_{i-\frac{1}{2}})\right)- 
h \, \bm{\mathcal{S}}_i=0,
%\end{equation}
\]
where $\bm{\mathcal{F}}(\phi_{i+\frac{1}{2},-},\phi_{i+\frac{1}{2},+};x_{i+\frac{1}{2}})$ is a two-point numerical flux evaluated at point
$x_{i+\frac{1}{2}}$. The scalar values $\phi_{i+\frac{1}{2},-}$ and $\phi_{i+\frac{1}{2},+}$ stand for the approximations
on the left and right sides of the interface $x_{i+\frac{1}{2}}$, respectively, and $\bm{\mathcal{S}}_i$ is an approximation 
of the mean value of the source term computed over cell $K_i$, for instance using a Gaussian quadrature formula.
Boundary conditions are taken into account in the first and last cells by the equations
$\phi_{\frac{1}{2},-}=\phi_\text{L}$ and $\phi_{I+\frac{1}{2},+}=\phi_\text{R}$.

\subsection{Cell polynomial degree and residual formulation}

Although polynomial reconstruction is an efficient tool to provide accurate approximations, one has to pay attention to the local regularity and choose the appropriate polynomial degree for the representation. 
For a cell $K_i$, we denote by $d_i$ its associated degree and $\mathcal M=(d_i)_{i=1,\ldots,I}$ stands for the Cell Polynomial Degree (CPD) map, that is, the vector which indicates the degree of the reconstruction for each cell. Moreover, we denote by 
\[
d_{i+\frac{1}{2}}=\min(d_i,d_{i+1})
\]
the effective polynomial degree we use to compute approximations on interface $x_{i+\frac{1}{2}}$.
A second important notion is the associated stencil $S_i$ for cell $K_i$ in relation with degree $d_i$. The simplest choice is the central stencil made of $\lceil d_i/2\rceil$ neighbor cells on each sides of $K_i$. We denote by $\mathcal{S}=(S_i)_{i=1,\ldots,I}$ the Cell Stencil (CS) map, that is the 
vector which indicates the reconstruction stencil for each cell.
% If the solution is regular enough according to the local sampling, we usually choose the central stencil $S_i^c$ made of $d_i/2$ neighbor cells on each sides of $K_i$. However in the case of an irregular solution or a steep gradient in the vicinity of $K_i$ (i.e within $S_i^c$), the polynomial reconstruction usually presents spurious numerical oscillations due to Gibbs phenomenon.  
Finally we set
\[
\phi_{i+\frac{1}{2},-}=\widehat{\phi}_{i}  (x_{i+\frac{1}{2}};d_{i+\frac{1}{2}}; S_i)\quad\text{and}\quad
\phi_{i+\frac{1}{2},+}=\widehat{\phi}_{i+1}(x_{i+\frac{1}{2}};d_{i+\frac{1}{2}}; S_{i+1}).
\] 
We introduce the residual $\mathcal{G}_i(\Phi, \mathcal M, \mathcal{S})$ at cell $K_i$ that is given by
\begin{eqnarray} \label{eq:residual_i}
\mathcal{G}_i(\Phi,\mathcal M, \mathcal S)=
\left(\bm{\mathcal{F}}(\phi_{i+\frac{1}{2},-},\phi_{i+\frac{1}{2},+};x_{i+\frac{1}{2}})-
\bm{\mathcal{F}}(\phi_{i-\frac{1}{2},-},\phi_{i-\frac{1}{2},+};x_{i-\frac{1}{2}})\right)
- h \, \bm{\mathcal{S}}_i,
\end{eqnarray}
and the global residual $\mathcal{G}(\Phi,\mathcal M, \mathcal S)=(\mathcal{G}_i(\Phi,\mathcal M, \mathcal S))_{i=1,\ldots,I}\in\mathbb R^I$ gathers expressions of the residuals of all cells. 

For given CPD and CS maps, $\mathcal M$, $\mathcal{S}$ respectively, we introduce the residual operator
\[
\Phi\in\mathbb R^I  \longrightarrow \mathcal{G}(\Phi,\mathcal M, \mathcal{S}) \in\mathbb R^I
\] 
and seek a solution $\Phi^{\mathcal M, \mathcal S}$ for the system of equations $\mathcal{G}(\Phi,\mathcal M, \mathcal{S})=0$ using an iterative  solver and given an initial guess. 

We also notice that the case where $d_i=0$ corresponds to the first-order approximation. On the other hand, if we are dealing with smooth solutions and $d_i=d$ with $d$ given by the user, then solving $\mathcal{G}(\Phi, \mathcal M, \mathcal{S})=0$ provides a $(d+1)$th-order approximation. However the point is to deal with solutions involving discontinuities. 
In this case with centered stencils, the CPD map admits high degrees in smooth regions but requires low ones where the solution is irregular.

Since we do not know \apriori the regularity of the solution, $\mathcal M$ also turns to be an unknown of the problem.
Accordingly, for smooth solutions, the stencil $S_i$ should be centered around cell $K_i$ while in the vicinity of discontinuous solution it is preferable to shift the stencil away. Doing so the reconstructed polynomial is less oscillatory and its associated linear/nonlinear system is also easier to solve.
Notice that the ENO/WENO family of schemes also use shifted stencils to build their polynomial reconstructions.
Therefore one looks for a solution constituted by the values of the solution itself but also the CPD and CS maps, $(\Phi^\star, \mathcal M^\star, \mathcal{S}^\star)$, such that $\mathcal{G}(\Phi^\star, \mathcal M^\star, \mathcal{S}^\star)=0$.

\subsection{Steady-state iterative solvers}

Regardless of the linearity of the PDE, 
%Even dealing with linear equations, 
the problem is nonlinear due to the choice of the CPD map $\mathcal M$. Therefore iterative procedures should produce a succession of approximations $(\Phi^k,\mathcal M^k, \mathcal{S}^k)$ that converges to a satisfying solution $(\Phi^\star,\mathcal M^\star, \mathcal{S}^\star)$. Furthermore, if the system of PDEs is nonlinear (B\"urgers or Euler equations), an extra inner loop is required to solve the nonlinear problem for fixed given CPD and CS maps $\mathcal M,\mathcal{S}$, respectively. 

We will consider two approaches to obtain the numerical approximations:
\begin{itemize}
\item \textit{Linear/Nonlinear solver (L/NL)}: in this approach we tackle directly $\mathcal{G}(\Phi, \mathcal M, \mathcal{S})=0$, con\-si\-de\-ring in this work a Gauss-type solver and OCTAVE nonlinear solver \texttt{fsolve} for the L and NL cases, respectively.

%Then new maps $\mathcal M^{k+1},\mathcal{S}^{k+1}$ are provided on the base of the data of $\Phi^{k+1}$ (the determination of these maps by MOOD is explained in the next section) and the iterative process proceeds up to convergence, {\it i.e.} $\mathcal M^k=\mathcal M^{k+1}$ and $\mathcal S^k=\mathcal S^{k+1}$ since in that case, we have $\Phi^k=\Phi^{k+1}$ and the steady-state approximate solution is obtained.
\item \textit{Time marching solver (TM)}: a second strategy consists in reformulating equation~\eqref{eq:generic_equation} by introducing a fictitious time dependency as
\begin{eqnarray}
\label{eq:linear_time}
\frac{\partial \widetilde \phi}{\partial t}+\frac{\partial\mathbf{F}(\widetilde \phi)}{\partial x} =\mathbf{S} 
\end{eqnarray}
with prescribed initial and boundary conditions. 
The time parameter $t$ is devoted to tend to infinity providing the steady state solution because the 
%searched 
steady-state solution is such that $\phi=\lim_{t \rightarrow \infty} \widetilde{\phi}(\cdot,t)$. In this work equation~\eqref{eq:linear_time} is solved using an explicit 
first order Runge-Kutta time discretization and high accurate space finite volume code as presented in~\cite{CDL3}.
% which computes the numerical FV solution $\widetilde \Phi^n=(\widetilde \phi_i^n)_{i=1,\ldots, I}$ where the maps $\mathcal M^n, \mathcal S^n$ %are adapted for each time step by the MOOD algorithm described in the next section.
%are fixed until two consecutive solutions do not differ anymore, that is
%\[
%\|\mathcal G(\tilde\Phi^{k,n+1})\|\leq\varepsilon
%\]
%\begin{eqnarray}
%\label{eq:linear_time2}
%\| \widetilde \phi^{n+1} - \widetilde \phi^{n} \| \leq \varepsilon,
%\end{eqnarray}
%with $\varepsilon$ of the order of the machine precision.
\end{itemize}

%
% END  GOVERNING EQUATIONS AND SCHEMES
%
%==================================================================================

%==================================================================================
%
% MOOD STENCIL
%
%\section{MOOD methodology} \label{sec:MOOD}
\section{The MOOD method} \label{sec:MOOD}

In~\cite{1D_SS} it is shown that the dynamical determination of the CPD map $\mathcal{M}$ is important to ensure stability, robustness, and accuracy simultaneously. 
In the case of a discontinuous solution, high accurate schemes produce numerical instabilities resulting into non-physical oscillations in the vicinity of the steep gradients.
But in the presence of a discontinuous solution or a steep gradient, the cell polynomial degree $d_i$ could drop to zero to avoid those spurious oscillations whereas the maximal user-defined degree $d^{\max}$ can be employed in the zones where the solution is regular enough, ensuring \textit{de facto} a local high accuracy of the approximate solution. 

The basis to dynamically determine the CPD map using a \aposteriori MOOD-like approach was developed in~\cite{CDL1,CDL2,CDL3}. 
The idea is the following: for a given stage $k$ and its associated map $\mathcal{M}^{k}$ with centered CS map $\mathcal{S}^{k}$, a candidate solution $\Phi^{k}$ is computed.
We then perform a detection procedure to identify which cells present numerical artifacts that would demand more dissipation. For those cells only we reduce the polynomial degree following a cascade, \textit{i.e.}, a decreasing sequence of polynomial degrees, maintaining the polynomial degrees of the good cells alike. This determines the updated CPD map $\mathcal{M}^{k+1}$ again with centered CS map $\mathcal{S}^{k}$ which is further employed to compute the new candidate solution $\Phi^{k+1}$.

In this section we briefly recall the three main ingredients of MOOD scheme from \cite{1D_SS}: detectors, cascade, and MOOD loop.

%% ---- FIG -----
%\begin{figure}[ht]
%	\begin{center} 
%		\includegraphics[width=0.6\textwidth]{figures/NLsolver1} 
%		\caption{ \label{fig:solvers}
%			Sketch of the original MOOD scheme from \cite{1D_SS}.
%			Starting from an initial guess $\Phi^0$ and CPD map $\mathcal{M}^0$ ($k=0$), the solver furnishes a candidate solution $\Phi^{k+1}$ which is tested (detection) and subsequent reduction (decrementing) of the polynomial degrees map may occur giving a new CPD map $\mathcal{M}^{k+1}$. If so, or if convergence is not reached then another iterate is performed.}
%	\end{center}
%\end{figure}
%% ---- FIG -----

\subsection{MOOD detectors} \label{ssec:detectors}

\begin{figure}
	\centering
	\includegraphics[width=1.0\textwidth]{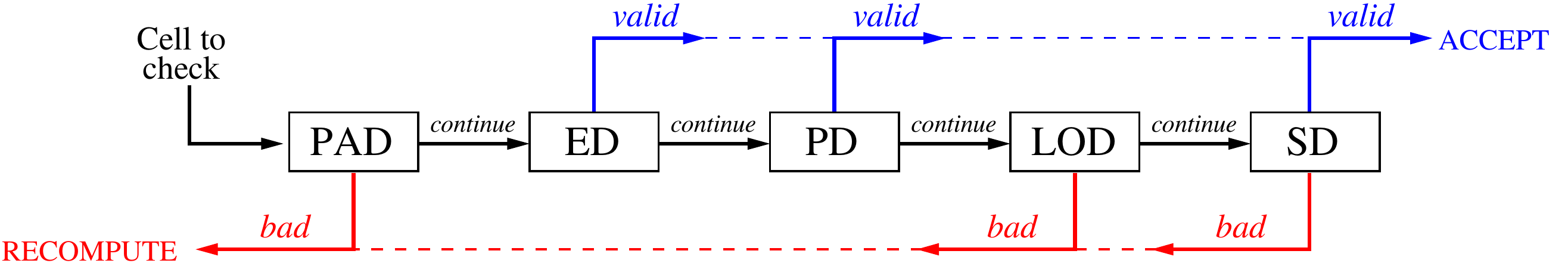}
	\caption{Chain of detectors used to check if a cell value is valid and accepted, or, if it is discarded 
		and needs to be recomputed after reduction of the local cell polynomial degree.
		Each of the top/bottom arrow permits to exit the chain prematurely.
		PAD: Physical Admissible Detector,
		ED: Extrema Detector,
		PD: Plateau Detector,
		LOD: Local Oscillation Detector,
		SD: Smoothness Detector.
	}
	\label{fig:chain_detector}
\end{figure}

The detection chain is the procedure by which the algorithm analyses the validity of a candidate solution in a given cell $K_i$.
If the cell is invalid, the solution needs to be locally recomputed with more dissipation,  by reducing the degree of the polynomial reconstruction associated to $K_i$. 
The detection chain is composed of a succession of elementary detector procedures, each focusing on a specific potential problem. 
We refer to \cite{1D_SS} for a complete description of the detectors and only describe them briefly for the sake of completeness, see Fig.~\ref{fig:chain_detector} for a sketch. 

% --- FIG --- chain detector
% --- FIG ---

\begin{description}
	\item[PAD] Physical Admissible Detector. The candidate solution satisfies the PAD condition on cell $K_i$ if it represents a relevant physical state 
	of the system of PDEs.
	Such a condition is crucial for the Euler system and it consists of the positivity of the density and pressure values.
	\item[ED] Extrema Detection. If a solution does not present a local extrema then the solution is accepted. 
	\item[PD] Plateau Detector. This criteria checks if the local curvature is too small (with respect to a tolerance $\varepsilon_\text{PD}=h$), if so the polynomial degree should be maintained as is.
	\item[LOD] Local Oscillation Detector. This criteria checks the variation of the signs of the local curvatures.
	\item[SD] Smoothness Detector. This detector determines if the minimum and the maximum curvatures are close enough with respect to the threshold parameter $\varepsilon_\text{SD}=0.25$, meaning that the nu\-me\-ri\-cal solution is considered as locally smooth enough.
\end{description}

\subsection{MOOD cascade} \label{ssec:MOOD_cascade}

The ``decrementing'' procedure is designed to reduce the polynomial degree $d_i$ when the cell has not passed the detectors.
To do so, we must determine a so-called ``cascade'' of polynomial degrees that the algorithm will successively try. A costly but systematic cascade consists in reducing the polynomial degree by one from $d^{\max}$ to zero. Notice that $d_i=d^{\max}$ leads to the most accurate scheme while $d_i=0$ to the most dissipative one.
A less expensive cascade may ignore some intermediate stages. 
For example in this work, following \cite{1D_SS}, we use $d^{\max}=5\rightarrow 2 \rightarrow 1 \rightarrow 0$ which ignores two intermediate stages but still has a nominal 6th order of accuracy on smooth flows.

\subsection{MOOD loop} \label{ssec:MOOD_loop}

The MOOD loop is the iterative procedure which, given a candidate solution $\Phi^{k+1}$ (obtained from a CPD and CS maps $\mathcal M^k,\mathcal S^k$, respectively) 
and an initial solution $\Phi^k$, states on its validity via the detection procedure and updates the CPD map via the decrementing procedure if the solution is invalid. 
This more appropriate CPD map $\mathcal M^{k+1}$ is further used along with initial solution $\Phi^{k+1}$ to compute a new candidate solution, which will be further tested for validity.

We present in Algorithm~\ref{L_NLsolver} the Linear/Nonlinear solvers and in Algorithms~\ref{TM1} and \ref{TM2} two versions of the Time Marching approach.

\begin{algorithm}
    \caption{\textit{L/NL solver}} \label{L_NLsolver}
    \KwData{$\Phi^0$, $\mathcal M^0$}
    \KwResult{$\Phi^\ast$, $\mathcal M^\ast$}

    $k\gets 0$\;
    \While{\text{TRUE}}{
    	$\Phi^{k+1}\gets \texttt{solver}(\Phi^k, \mathcal M^k)$\;
    	$\mathcal M^{k+1}\gets \texttt{detect\_and\_decrement}(\Phi^{k+1}, \mathcal M^k)$\;
		\If{$\mathcal M^{k+1}=\mathcal M^k$}{
			$\Phi^\ast\gets \Phi^{k+1}$, $\mathcal M^\ast\gets\mathcal M^k$\;
			\Return{$\Phi^\ast$, $\mathcal M^\ast$}\;
		}
		\Else{
			$k\gets k+1$\;
		}
    }
\end{algorithm}

\begin{algorithm}
    \caption{\textit{TM1 solver}} \label{TM1}
    \KwData{$\tilde\Phi^0$, $\mathcal M^0$, $\varepsilon$}
    \KwResult{$\Phi^\ast$, $\mathcal M^\ast$}

    $n\gets 0$\;
    \While{\text{TRUE}}{
	    $i\gets 0$\;
	    $\tilde\Phi^{n,0}\gets\tilde\Phi^n$, $\mathcal M^{n,0}\gets\mathcal M^0$\;
	    \While{\text{TRUE}}{
	    	$\tilde\Phi^{n,i+1}\gets \texttt{rk1}(\tilde\Phi^{n,i}, \mathcal M^{n,i})$\;
	    	$\mathcal M^{n,i+1}\gets \texttt{detect\_and\_decrement}(\tilde\Phi^{n,i+1}, \mathcal M^{n,i})$\;
			\If{$\mathcal M^{n,i+1}=\mathcal M^{n,i}$}{
				$\tilde\Phi^{n+1}\gets \tilde\Phi^{n,i+1}$\;
			}
			\Else{
				$i\gets i+1$\;
			}
	    }
	    \If{$\|\mathcal G(\tilde\Phi^{n+1})\|\leq\varepsilon$}{
		    $\Phi^\ast\gets \tilde\Phi^{n+1}$, $\mathcal M^\ast\gets\mathcal M^{n,i+1}$\;
		    \Return{$\Phi^\ast$, $\mathcal M^\ast$}\;
	    }
	    \Else{
		    $n\gets n+1$\;
	    }
    }
\end{algorithm}

\begin{algorithm}
    \caption{\textit{TM2 solver}} \label{TM2}
    \KwData{$\tilde\Phi^0$, $\mathcal M^0$, $\varepsilon$}
    \KwResult{$\Phi^\ast$, $\mathcal M^\ast$}

    $k\gets 0$\;
    \While{\text{TRUE}}{
        
	    $n\gets 0$\;
	    $\tilde\Phi^{k,n}\gets\tilde\Phi^k$\;
	    \While{\text{TRUE}}{
	    	$\tilde\Phi^{k,n+1}\gets \texttt{rk1}(\tilde\Phi^{k,n}, \mathcal M^{k})$\;
			\If{$\|\mathcal G(\tilde\Phi^{k,n+1})\|\leq\varepsilon$}{
				$\Phi^{k+1}\gets \tilde\Phi^{k,n+1}$\;
			}
			\Else{
				$n\gets n+1$\;
			}

	    }
    	
    	$\mathcal M^{k+1}\gets \texttt{detect\_and\_decrement}(\Phi^{k+1}, \mathcal M^k)$\;
		\If{$\mathcal M^{k+1}=\mathcal M^k$}{
			$\Phi^\ast\gets \Phi^{k+1}$, $\mathcal M^\ast\gets\mathcal M^k$\;
			\Return{$\Phi^\ast$, $\mathcal M^\ast$}\;
		}
		\Else{
			$k\gets k+1$\;
		}
    }
\end{algorithm}

% END   MOOD STENCIL
%
%==================================================================================

%==================================================================================
%
% DYNAMICAL STENCIL
%
%\section{Dynamical stencil methodology} \label{sec:as}
%
% CORE OF THE PAPER 
%
\section{The MOOD method supplemented with Adaptive Stencil} \label{sec:as}

In this section we illustrate why an Adaptive Stencil (AS) procedure during the polynomial reconstruction may be useful to increase the accuracy of the MOOD scheme. Next we propose an AS algorithm to be associated to the polynomial degree selection of the MOOD scheme.

\subsection{Goal and idea} \label{ssec:goal}

\begin{figure}
	\centering
	\includegraphics[width=0.329\textwidth]{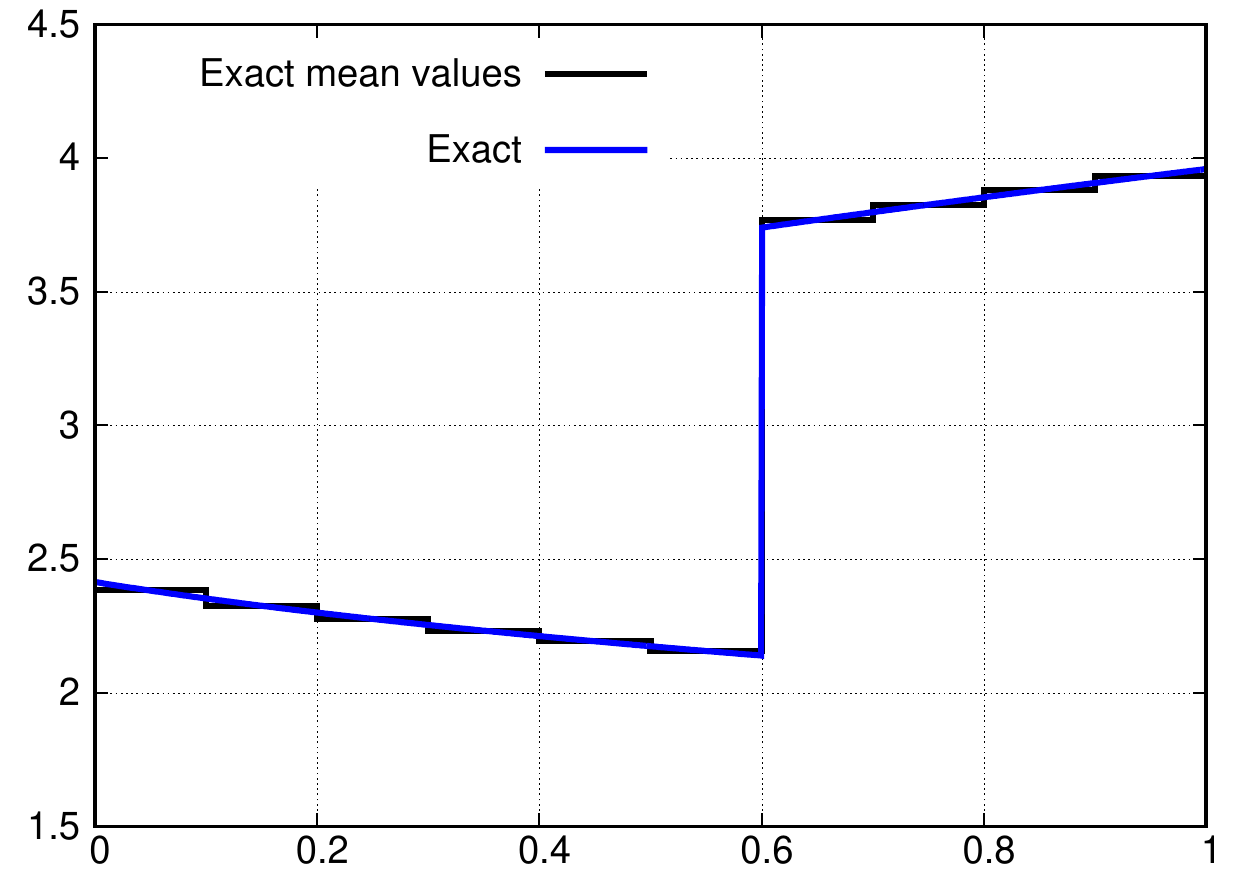}
	\includegraphics[width=0.329\textwidth]{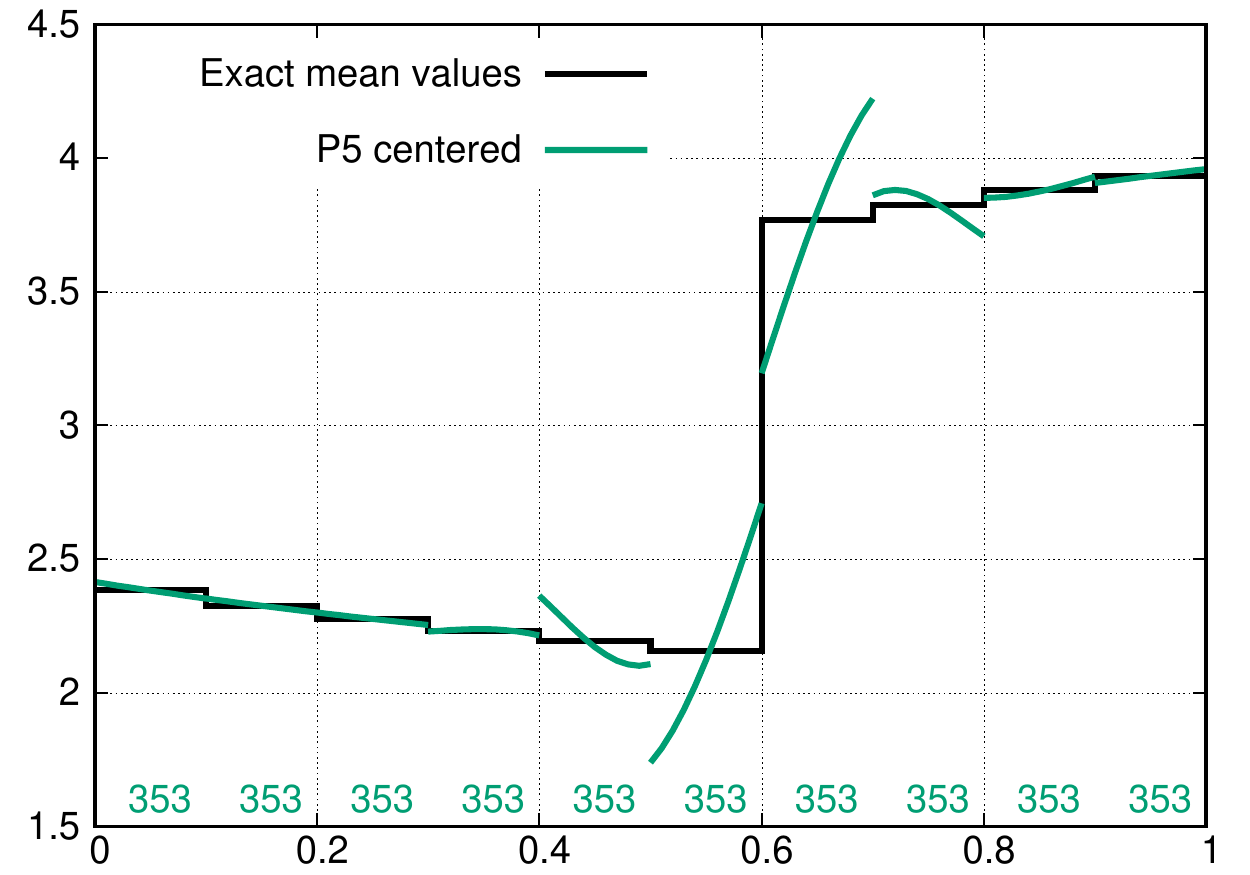}
	\includegraphics[width=0.329\textwidth]{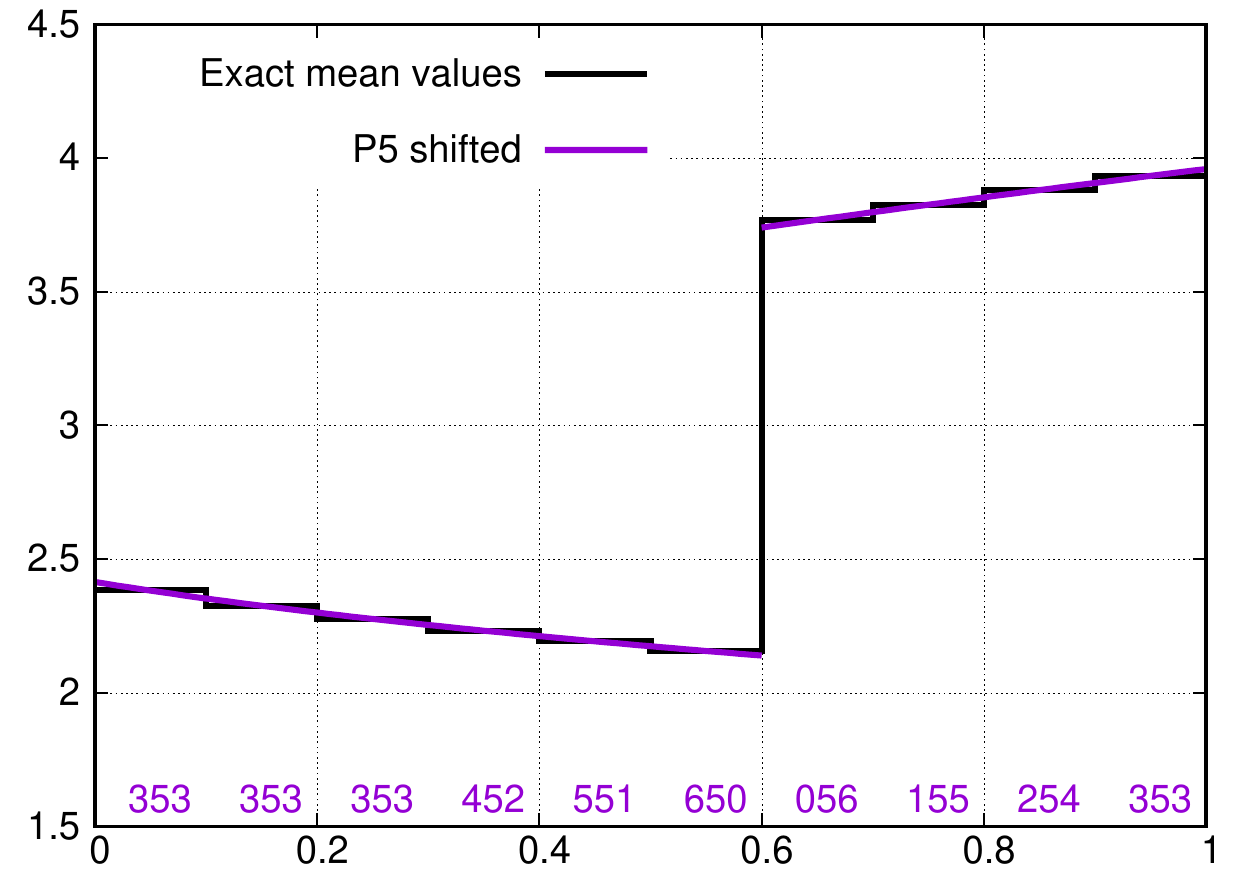}
	\caption{
		Left panel: exact function (blue) composed of two branches separated by a discontinuity along with its associated mean values on a 
		$10$ cell uniform mesh ---
		Middle panel: centered stencil $\mathbb{P}_5$ reconstructions, where label $ldr=353$ refers to $l=r=3$ neighbor cells to the left and right of the current cells for a $d=5$ polynomial degree ---
		Right panel: shifted stencil $\mathbb{P}_5$ reconstruction, the stencils are chosen such that the discontinuity is never included into it.
	}
	\label{fig:DSexample}
\end{figure}

The MOOD scheme described above and in \cite{1D_SS} employs only centered stencils. As such only the CPD map is modified according to the numerical solution via the detection procedure. Indeed the polynomial degree is decremented (and a centered stencil is always employed) until the solution is well behaved. The main drawback of this strategy is the exclusive usage of centered stencil whereas, for some situations, maintaining the maximal polynomial degree but shifting the stencil would be sufficient to ensure the validity of the candidate solution.
% Illustrative example
As an illustrative example in Fig.~\ref{fig:DSexample} we present the centered or shifted stencil $\mathbb{P}_5$ reconstructions of an exact function composed of two non-polynomial branches and a discontinuity. The left panel shows the exact solution (blue) and the associated $10$ mean values (black) on a uniform mesh from which we obtain the two reconstructions; one using a centered stencil (middle panel in green) and the second one using a shifted stencil (right panel in purple). 
The centered $\mathbb{P}_5$ stencil reconstructions present some over/under-shoots in the vicinity of the discontinuity. Indeed, given a cell $i$, the centered stencil is formed by $l=3$ neighbor cells on the left and $r=3$ on the right; This is the meaning of labels $ldr=353$. As such the discontinuity located at the interface between the $6$th and $7$th cells does influence the polynomial reconstructions from the $4$th up to the $9$th cells. 
Contrarily, when some appropriate shifted stencils are employed, it is possible to reconstruct $\mathbb{P}_5$ polynomials without any spurious effect generated by the presence of the discontinuity. Notice that the stencils (and labels in the figure) are shifted to the left for cell $4$, $5$, and $6$ and to the right for cells $7$, $8$, and $9$. All other cells use their centered stencils.
% --- FIG --- DS example
% --- FIG ---

In this work we supplement the scheme with the possibility of adapting the CPD map and the reconstruction stencil accordingly. Let us start with the algorithm used to determine an adapted stencil given a CPD map.

\subsection{Adaptive cell stencil} \label{ssec:AS_algo}

Given a valid CPD map $\mathcal M$ in 1D, we can design an algorithm which determines the best shifted stencil to reconstructed $\mathbb{P}_d$ polynomials (if possible), $d$ being the maximal polynomial degree. Starting from a cell index $i$ and its CPD $d_i$, which may range from $0$ to $d$ the maximal possible CPD, the algorithm should find at least $d$ neighbor cells. We present in Algorithm~\ref{as_alg} a simple proposal to compute $\mathcal S^{k+1}_i$ with $\#\mathcal S^{k+1}_i=N$. This algorithm produces the stencils depicted in Fig.~\ref{fig:DSexample} (right panel).

\begin{algorithm}
    \caption{ \label{algo:S} Calculate $\mathcal S^{k+1}_i$} \label{as_alg}
    \KwData{$i$, $I$, $N$, $\mathcal M^{k+1}$}
    \KwResult{$\mathcal S^{k+1}_i$}

	$\mathcal S^{k+1}_i \gets \{i \}$\;
	$i_\text{L} \gets i-1$, $c_\text{L} \gets 0$\;
	$i_\text{R} \gets i+1$, $c_\text{R} \gets 0$\;

	\While{$\#\mathcal S^{k+1}_i<N$}{
		\If{$i_\text{L}=0$}{
			$\mathcal S^{k+1}_i \gets \mathcal S^{k+1}_i \cup \{\text{right cells indices}\}$; break\;
		}

		\If{$i_\text{R}=I+1$}{
			$\mathcal S^{k+1}_i \gets \mathcal S^{k+1}_i \cup \{\text{left cells indices}\}$; break\;
		}

		\If{$((\mathcal M^{k+1}(i_\text{L})>\mathcal M^{k+1}(i_\text{R})) \lor (\mathcal M^{k+1}(i_\text{L})=\mathcal M^{k+1}(i_\text{R}) \land c_\text{L} \leq c_\text{R})$}{
			$\mathcal S^{k+1}_i \gets \mathcal S^{k+1}_i \cup \{i_\text{L}\}$\;
			$i_\text{L} \gets i_\text{L}-1$, $c_\text{L} \gets c_\text{L}+1$\;
		}
		\Else{
			$\mathcal S^{k+1}_i \gets \mathcal S^{k+1}_i \cup \{i_\text{R}\}$\;
			$i_\text{R} \gets i_\text{R}-1$, $c_\text{R} \gets c_\text{R}+1$\;
		}
	}
	remove $i$ from $\mathcal S^{k+1}_i$
\end{algorithm}

The principle is based on the CPD map which plays the role of a smoothness indicator. Given a CPD map $\mathcal M^{k+1}$ over the grid, we produce stencil $\mathcal S_i^{k+1}$ for cell $i$ of prescribed size $N$. We first start with the singleton $\{i\}$, the index of the reference cell. Then we choose the adjacent cell (between the left and right cells) that have the largest degree. A priority criteria is considered if the two cells have the same degree. We obtain a two-cell stencil and start again the search with the adjacent cells. At each iteration, we agglomerate the best adjacent cell to the current stencil and we stop the procedure when we reach the prescribed size $N$. When the stencil is complete, we remove the reference cell $i$ to fulfill the definition given in section \ref{sec::def_stencil}. The algorithm shifts the stencil in direction away from any discontinuity and merge cells from the regular region.  

\subsection{MOOD+AS algorithm} \label{ssec:MOOD_AS}

Regardless of the solver used in the MOOD solver (non-linear or time-marching one), the MOOD+AS algorithm is given by the following procedure. A solution, denoted as $\Phi^{\text{MOOD}}$, is computed with the original MOOD method  \cite{1D_SS} using centered stencils. The solution also provides the CPD map $\mathcal{M}^{\text{MOOD}}$, that characterises the regularity. Such a solution is acceptable by construction since the CPD map is computed to pass the detection criteria.

We then carry out the successive operations 
\[
(\Phi^{k},\mathcal{M}^{d^{\max}},\mathcal{S}^{k})\xrightarrow[]{MOOD}
(\Phi^{k+1},\mathcal{M}^{k+1})\xrightarrow[]{AS}
\mathcal{S}^{k+1}
\]
where we initialize the loop with $\Phi^{0}=\Phi^{\text{MOOD}}$ and $\mathcal{M}^{0}=\mathcal{M}^{\text{MOOD}}$. The first stage corresponds to the MOOD solver using the collection of stencils $\mathcal S^{k}$ and produce an admissible solution $\Phi^{k+1}$ with its CPD map $\mathcal{M}^{k+1}$. Then, we apply Algorithm~\ref{algo:S} to adapt the stencils according to the CPD map $\mathcal{M}^{k+1}$ and get the new collection of stencils $\mathcal{S}^{k+1}$. The loop stops when two successive CDP map are equal, namely $\mathcal{M}^{k}=\mathcal{M}^{k+1}$. It is worth mentioning that the two operators are independent and  implemented in separeted routines. Hence, the AS strategy is quite versatile and could be used in an other context.

% END   DYNAMICAL STENCIL
%
%==================================================================================

%==================================================================================
%
% NUMERICS
%
%\section{Numerical results} \label{sec:numerics}
\section{Numerical results} \label{sec:numerics}
% Models
Numerical experiments are carried out with the MOOD+AS scheme to give evidences of the method efficiency, namely to provide a better accuracy still preserving the robustness. We detail three relevant problems that cover a large spectrum of applications. The linear advection equation with a regular and a non-differential velocity is a classical building block in numerous problems. Next we analyse the method for the well-known B\"urgers equation involving a simple non-linearity. At last, we deal with the Euler system that represents the standard vector-valued non-linear hyperbolic system.

% MOOD scheme 6th order
As already mentioned we only focus on one single cascade $d^{\max}=5\rightarrow 2 \rightarrow 1 \rightarrow 0$ ignoring two intermediate stages and having a nominal 6th order of accuracy on smooth flows.
% Diagnostics
To assess the error and convergence rate, we use the $L^1$ and $L^\infty$ errors computed by
\[
E_1 \equiv E_1(\Phi,I)=\Sum_{i=1}^{I} |\phi_i-\phi_i^{\text{ex}}| h,
\qquad
E_{\infty}\equiv E_\infty(\Phi,I) =\displaystyle\max_{i=1}^{I} |\phi_i-\phi_i^{\text{ex}}|.
\]
We further define the rate of convergence between two solutions/meshes $(\Phi_k,I_k)$, for $k=1,2$ where $I_1<I_2$ as
\[
\mathcal{O}_{\alpha} \equiv \mathcal{O}_{\alpha} \left( (\Phi_1,I_1);(\Phi_2,I_2) \right) = 
\frac{|\log E_\alpha(\Phi_1,I_1)/E_\alpha(\Phi_2,I_2)|}{|\log I_1/I_2|}, \qquad
\alpha=1,\infty.
\]

The adaptive stencil strategy implies that one has to characterize the cells' position belonging to the stencil regarded from the reference cell where the polynomial reconstruction is evaluated. Table~\ref{colors_notation} gives the color codes used in the CPD map representations to distinguish 
%centered and shifted 
different stencils. For example, a centered stencil (same number of white circles cells on the left and right side of the reference cell marked by a coloured circle) is colored in green while a full left upwind stencil (all the cells are on the left side) is colored in violet. Obviously the colors' code also depends on the polynomial degree of the reconstruction.

\begin{table}\centering
\caption{Colors notation for stencils.}
\label{colors_notation}
\begin{tabular}{@{}ccccc@{}}\toprule
$d$ & $\# S_i$ & $l$ & $r$ & \\\midrule
1 & 2 & 2 & 0 & \CPDwhite\CPDwhite\CPDblue\CPDwhite\textcolor{white}{O}\\
   &   & 1 & 1 & \textcolor{white}{O}\CPDwhite\CPDgreen \CPDwhite\textcolor{white}{O}\\
   &   & 0 & 2 & \textcolor{white}{O}\textcolor{white}{O}\CPDyellow \CPDwhite\CPDwhite\\\midrule
2 & 2 & 2 & 0 & \CPDwhite\CPDwhite\CPDblue\textcolor{white}{O}\textcolor{white}{O}\\
   &   & 1 & 1 & \textcolor{white}{O}\CPDwhite\CPDgreen \CPDwhite\textcolor{white}{O}\\
   &   & 0 & 2 & \textcolor{white}{O}\textcolor{white}{O}\CPDyellow \CPDwhite\CPDwhite\\\midrule
5 & 6 & 6 & 0 & \CPDwhite\CPDwhite\CPDwhite\CPDwhite\CPDwhite\CPDwhite\CPDviolet \textcolor{white}{O}\textcolor{white}{O}\textcolor{white}{O}\textcolor{white}{O}\textcolor{white}{O}\textcolor{white}{O}\\
   &   & 5 & 1 & \textcolor{white}{O}\CPDwhite\CPDwhite\CPDwhite\CPDwhite\CPDwhite\CPDindigo \CPDwhite\textcolor{white}{O}\textcolor{white}{O}\textcolor{white}{O}\textcolor{white}{O}\textcolor{white}{O}\\
   &   & 4 & 2 & \textcolor{white}{O}\textcolor{white}{O}\CPDwhite\CPDwhite\CPDwhite\CPDwhite\CPDblue \CPDwhite\CPDwhite\textcolor{white}{O}\textcolor{white}{O}\textcolor{white}{O}\textcolor{white}{O}\\
   &   & 3 & 3 & \textcolor{white}{O}\textcolor{white}{O}\textcolor{white}{O}\CPDwhite\CPDwhite\CPDwhite\CPDgreen \CPDwhite\CPDwhite\CPDwhite\textcolor{white}{O}\textcolor{white}{O}\textcolor{white}{O}\\
   &   & 2 & 4 & \textcolor{white}{O}\textcolor{white}{O}\textcolor{white}{O}\textcolor{white}{O}\CPDwhite\CPDwhite\CPDyellow \CPDwhite\CPDwhite\CPDwhite\CPDwhite\textcolor{white}{O}\textcolor{white}{O}\\
   &   & 1 & 5 & \textcolor{white}{O}\textcolor{white}{O}\textcolor{white}{O}\textcolor{white}{O}\textcolor{white}{O}\CPDwhite\CPDorange \CPDwhite\CPDwhite\CPDwhite\CPDwhite\CPDwhite\textcolor{white}{O}\\
   &   & 0 & 6 & \textcolor{white}{O}\textcolor{white}{O}\textcolor{white}{O}\textcolor{white}{O}\textcolor{white}{O}\textcolor{white}{O}\CPDred \CPDwhite\CPDwhite\CPDwhite\CPDwhite\CPDwhite\CPDwhite\\\bottomrule
\end{tabular}
\end{table}

\subsection{Advection equation} \label{ssec:advection}

The scalar linear convective steady-state problem reads
\[
\frac{\textrm{d}}{\textrm{d}x}\Big (u(x)\phi(x)\Big)=S(x),\quad x\in \Omega,
\]
where $u\equiv u(x)>0$ is a given velocity function and Dirichlet boundary condition $\phi_{\text{L}}$ is prescribed on the left side of domain $\Omega$. Since $u$ may present some irregularities, the solution is not necessarily smooth enough to perform all the polynomial reconstructions with $d^{\max}$ hence the CPD and CS maps must be adapted accordingly.

The numerical flux $\mathcal F(\phi_{i+\frac{1}{2},-},\phi_{i+\frac{1}{2},+};x_{i+\frac{1}{2}})$ is the classical upwind one.
For given CPD and CS maps $\mathcal M, \mathcal S$, and noticing that the problem $\mathcal{G}(\Phi, \mathcal M, \mathcal S)=0$ is affine, we use a Gauss-type procedure to nullify the residual.
% Chain detector and cascade
Next, the chain detector associated to this linear equation follows the chain: ED, PD, LOD, and SD where all these detectors have been defined in the previous section (see also \cite{1D_SS}). 

%--- F I G ------
\begin{figure}\centering
\begin{tabular}{@{}cc@{}}
\includegraphics[width=0.4\textwidth]{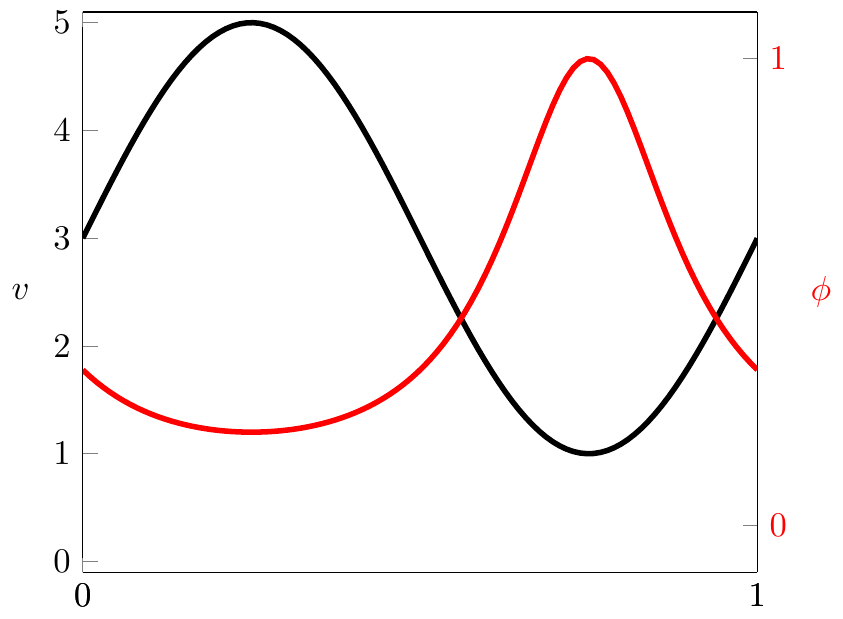} &
\includegraphics[width=0.42\textwidth]{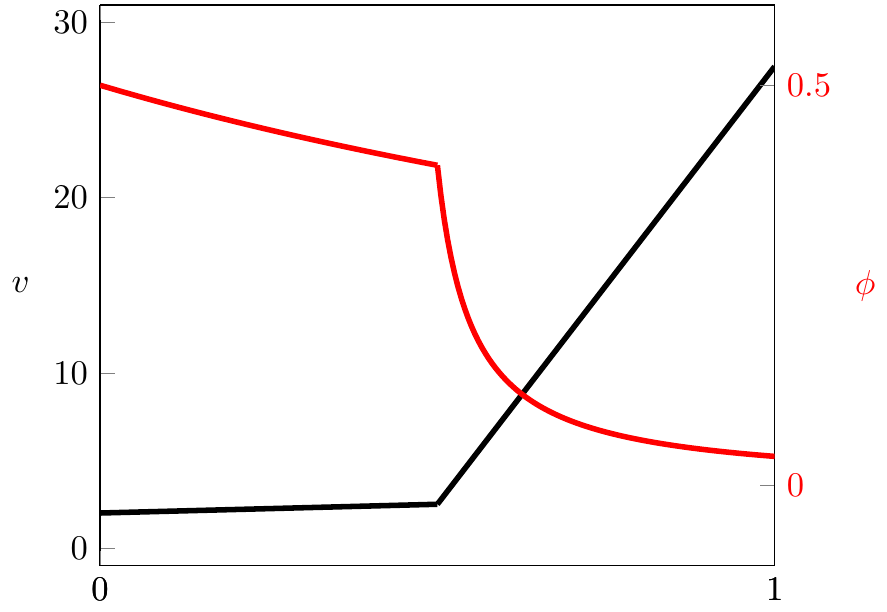}
\end{tabular}
\caption[]{\label{fig:ex_sol_linear} Exact solutions (\phiphi) 
and associate velocities (\velocity) for the linear equation on $\overline\Omega=[0;1]$ --- 
$\phi(x)=1/u(x)$ ---
Left: regular solution $u(x)= 2 \sin(2\pi x) + 3$ --- 
Right: irregular solution $u(x)=x+2, \; \text{if} \;  x<1/2$ or $u(x)=50(x-1/2)+5/2$, otherwise.}
\end{figure}
%--- F I G ------

\subsubsection{Regular solution}
\begin{figure}\centering
\begin{tabular}{@{}c@{}c@{}c@{}}\toprule
exact and numerical solutions & CPD map & errors\\\midrule
\includegraphics[width=0.3\textwidth,align=c]{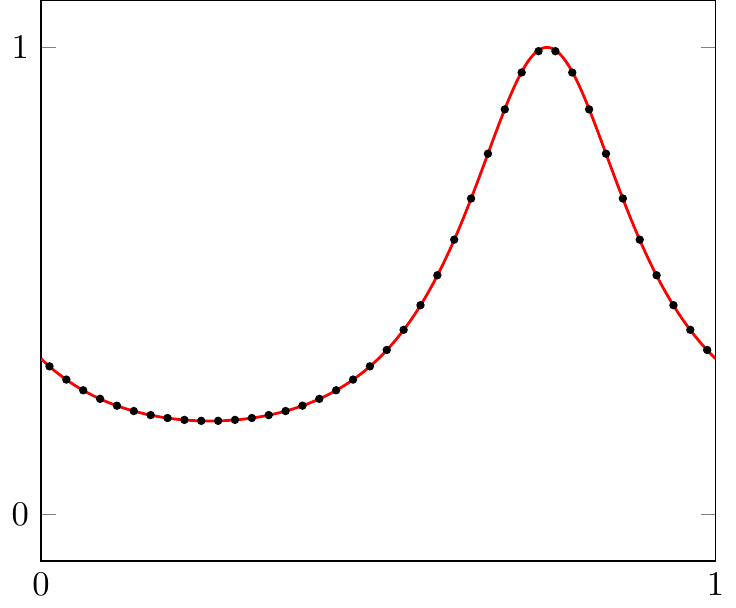} &
\includegraphics[width=0.3\textwidth,align=c]{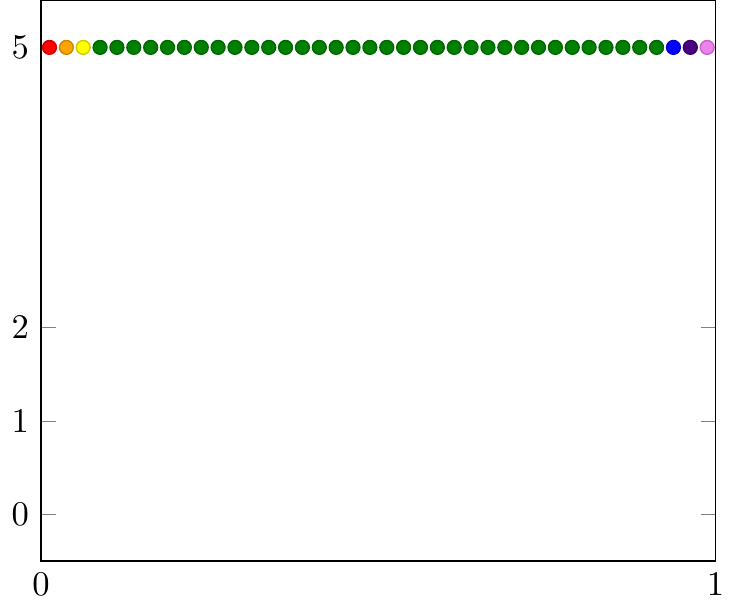} &
\includegraphics[width=0.3\textwidth,align=c]{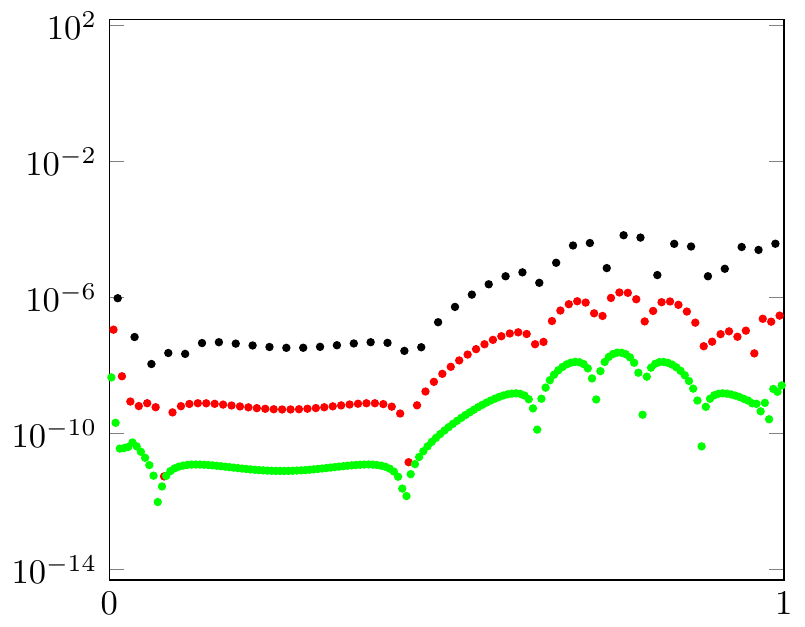}\\\bottomrule
\end{tabular}
\caption[]{\label{fig:L_ex1_L} Advection example 1 (regular): L solver with MOOD --- exact and numerical solutions for $I=40$ (left), CPD map for $I=40$ (center), and cell errors for $I=40$ (\blackblack), $I=80$ (\redred), and $I=160$ (\greengreen) (right).}
\end{figure}

A regular situation is achieved with the smooth velocity function $u(x)= 2 \sin(2\pi x) + 3$ leading to the solution depicted in Fig.~\ref{fig:ex_sol_linear} on left panel. The boundary conditions are set to $1/3$ at $x=0$ and an outflow condition at $x=1$.
In Table~\ref{tab:conv_ex1} one reports the $L_1$ and  $L_\infty$ errors, together with the convergence rates for a 6th order MOOD scheme. Optimal orders is reported for the L solver while we observe a small discrepancy with the TM1 and TM2 methods. We only report the sixth-order scheme results since similar evidences are obtained for the second- and fourth-order cases. 
Moreover Fig.~\ref{fig:L_ex1_L} presents the numerical solutions (left panel), the CPD maps (middle panel), and cell errors (right panel) for L scheme using the 6th-order MOOD algorithm. We observe that the Cell Polynomial Degree is maintained at its maximum leading to the optimal accuracy. No oscillations are detected. 
We do observe that the overall errors drops by almost two orders of magnitude from one mesh to the next refined one.
We do not present the TM1 and TM2 figures since they are almost similar to the L one.

\begin{table}\centering
\setlength\tabcolsep{3 pt}
{\tablepolice
\caption{\label{tab:conv_ex1} Advection with regular velocity: errors and convergence rate: L solver (left), TM1 solver (middle), and TM2 solver (right).}
\begin{tabular}{@{}cccccccc@{}}\toprule
$I$ & $E_1$ & ${\mathcal O}_1$ & $E_\infty$ & ${\mathcal O}_\infty$\\\midrule
40 & 1.1E$-$05 & --- & 6.9E$-$05 & --- \\
80 & 1.7E$-$07 & 6.0 & 1.4E$-$06 & 5.6 \\
160 & 2.6E$-$09 & 6.0 & 2.4E$-$08 & 5.9\\\bottomrule
\end{tabular}
\hskip 2em
\begin{tabular}{@{}cccccccc@{}}\toprule
$I$ & $E_1$ & ${\mathcal O}_1$ & $E_\infty$ & ${\mathcal O}_\infty$\\\midrule
40 & 1.1E$-$05 & --- & 6.6E$-$05 & --- \\
80 & 2.3E$-$07 & 5.5 & 1.8E$-$06 & 5.2 \\
160 & 6.2E$-$09 & 5.2 & 3.8E$-$08 & 5.5 \\\bottomrule
\end{tabular}
\hskip 2em
\begin{tabular}{@{}cccccccc@{}}\toprule
$I$ & $E_1$ & ${\mathcal O}_1$ & $E_\infty$ & ${\mathcal O}_\infty$\\\midrule
40 & 1.1E$-$05 & --- & 6.6E$-$05 & --- \\
80 & 2.3E$-$07 & 5.5 & 1.8E$-$06 & 5.2 \\
160 & 6.2E$-$09 & 5.2 & 3.8E$-$08 & 5.5 \\\bottomrule
\end{tabular}
}
\end{table}

\subsubsection{Irregular solution}
The second test is dedicated to the situation generated by an irregular velocity
\[
u(x) =
\begin{cases} 
x+2,           & \text{if $x\in[0;1/2]$},\\
50(x-1/2)+5/2, & \text{if $x\in[1/2;1]$},
\end{cases} 
\]
with the boundary conditions set to $1/2$ at $x=0$ and an outflow condition at $x=1$.
It leads to the solution depicted in Fig.~\ref{fig:ex_sol_linear} on right panel.
This solution is continuous but not differentiable at location $x=1/2$.

\begin{figure}\centering
\begin{tabular}{@{}cc@{}c@{}c@{}}\toprule
& exact and numerical solutions & CPD map & errors\\\midrule
MOOD & \includegraphics[width=0.3\textwidth,align=c]{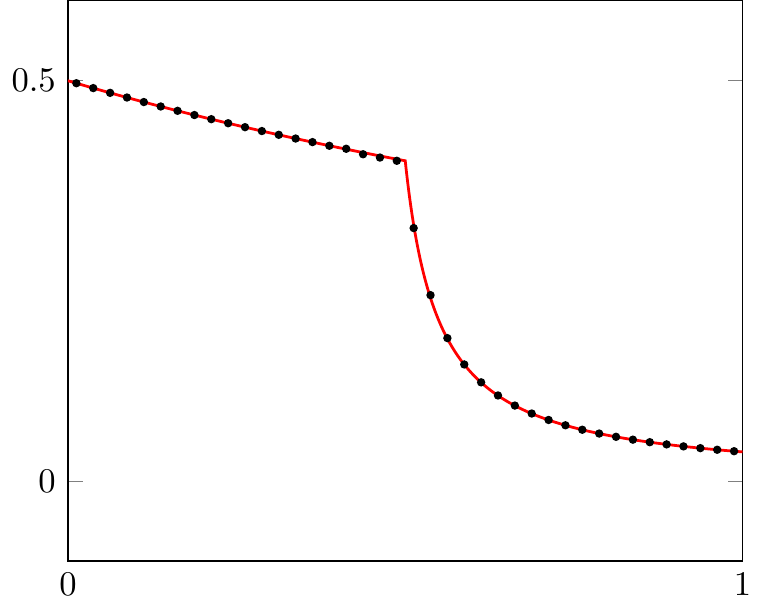} &
\includegraphics[width=0.3\textwidth,align=c]{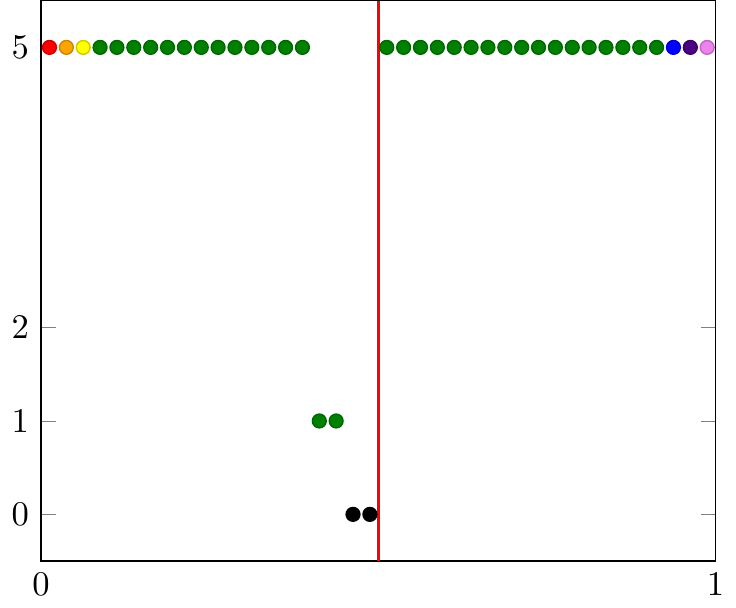} &
\includegraphics[width=0.3\textwidth,align=c]{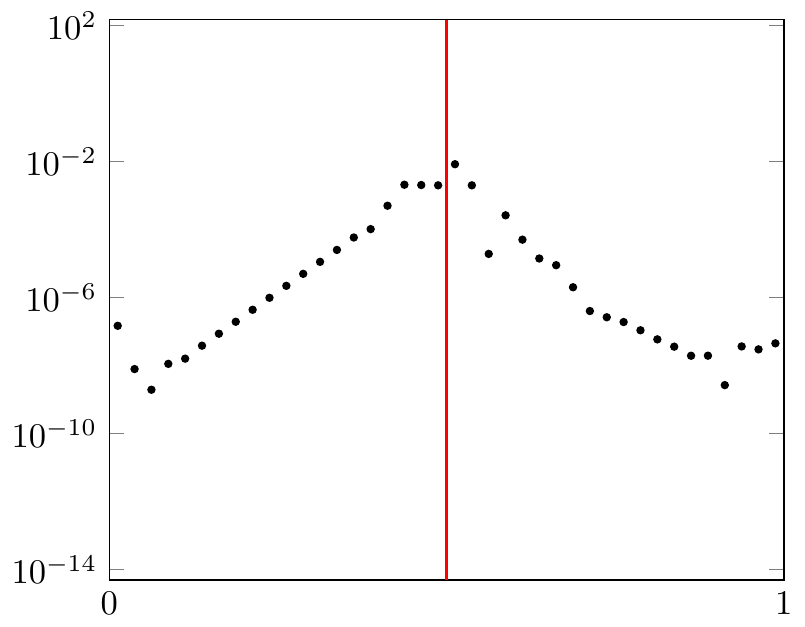}\\
AS & \includegraphics[width=0.3\textwidth,align=c]{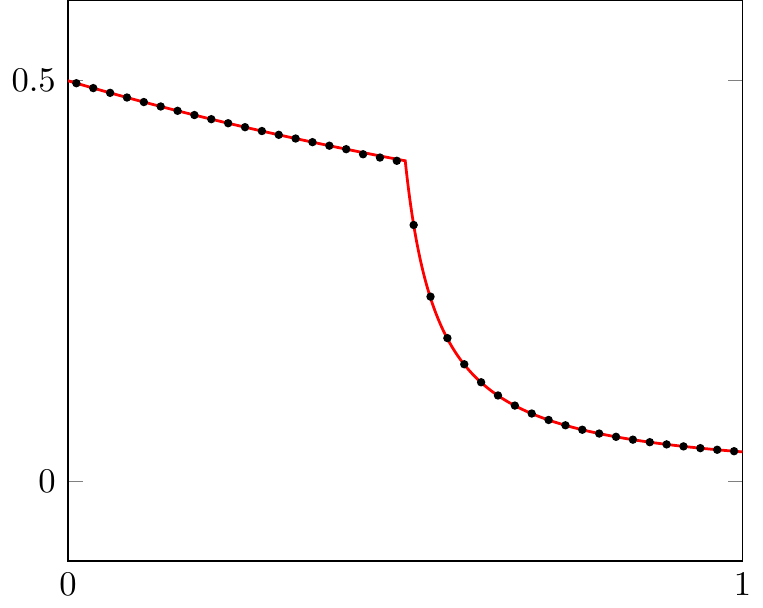} &
\includegraphics[width=0.3\textwidth,align=c]{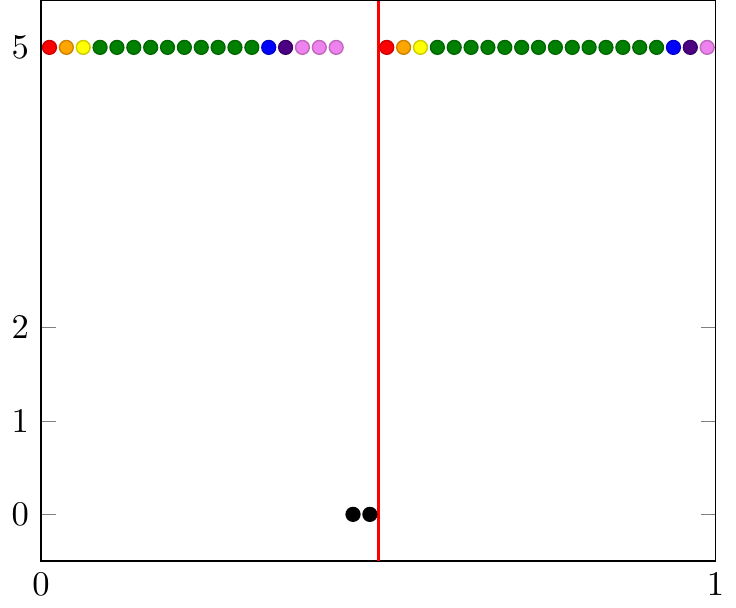} &
\includegraphics[width=0.3\textwidth,align=c]{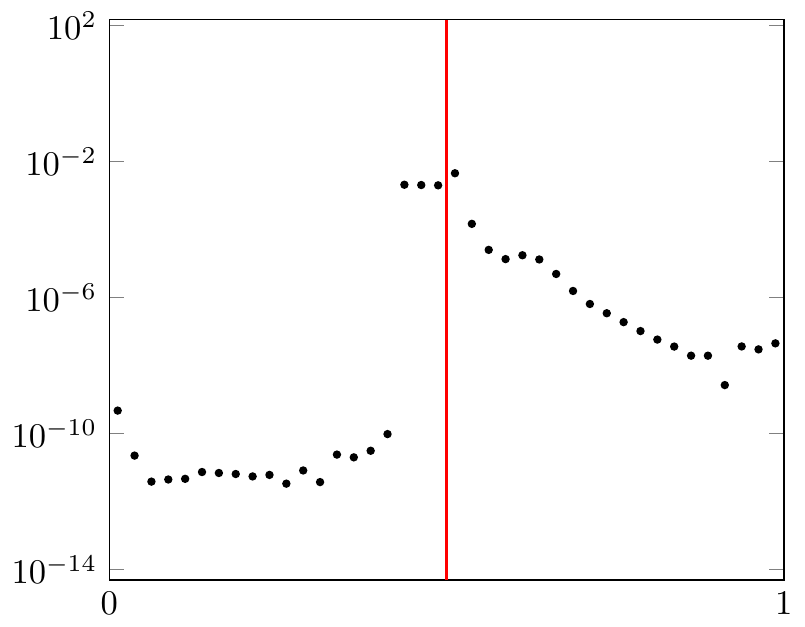}\\\bottomrule
\end{tabular}
\caption{\label{fig:L_ex2_L} Advection example 2 (irregular): L solver for $I=40$.}
\end{figure}

Unlimited reconstructions, \textit{i.e.}, using the systematic $\mathbb P_5$ polynomial representation, produces a non-physical solution and spurious oscillations are generated. To remove them, 
%overshooting, 
the original MOOD \textit{a posteriori} strategy (Static Stencil) is applied by reducing the polynomial degree for the pro\-ble\-ma\-tic cells. 
%that require a cure. 
We recall that the stencils are centered with respect to the reference cell
and it results in decrementing 
% leading to a strong cut 
of the polynomial degrees in the vicinity of 
a discontinuity.
%a shock. The alternative 
The proposed Adaptative Stencil method
% we propose 
contrarily modifies the stencil structure in order to select the neighbor cells %that enjoy 
compatible with a high degree of reconstruction 
%that helps 
helping to maintain the optimal accuracy.  

We present in Fig.~\ref{fig:L_ex2_L} the two strategies and their consequences for the accuracy and stability.
The first row displays the numerical solution, the CPD map, and the cell errors using the Static Stencil (the original MOOD method) whereas the second row 
%displays the same figures for 
presents the Adaptative Stencil results.
%strategies. 
In both cases, the stability is achieved and no oscillations are reported. The CPD map configurations are different since the Adaptative Stencil only suffers from two degree reduction cells. We notice the stencil upwinding close to the singularity in accordance with the relative position (left upwinded on the left side and accordingly on the right). The last column depicts the error values, cell by cell, and shows the 
%dramatic 
clear improvement. 
%of the accuracy.     

\begin{figure}\centering
\begin{tabular}{@{}cc@{}c@{}c@{}}\toprule
& exact and numerical solutions & CPD map & errors\\\midrule
MOOD & \includegraphics[width=0.3\textwidth,align=c]{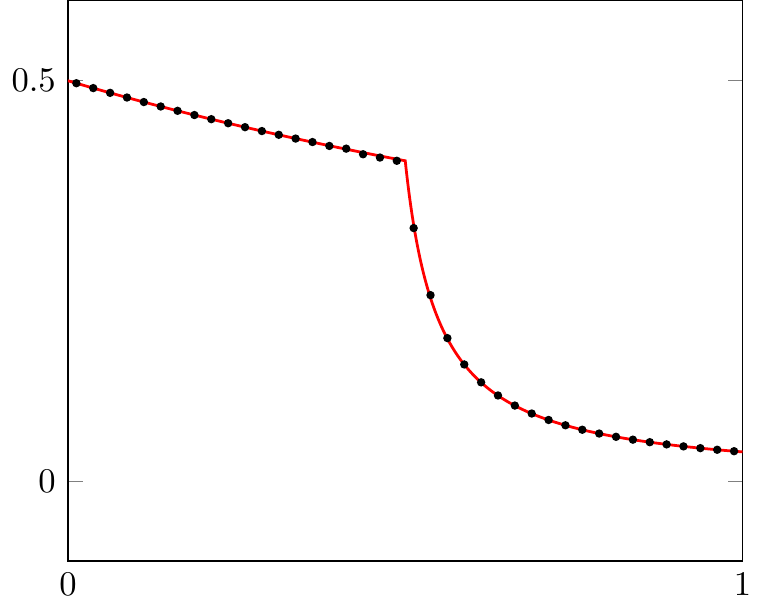} &
\includegraphics[width=0.3\textwidth,align=c]{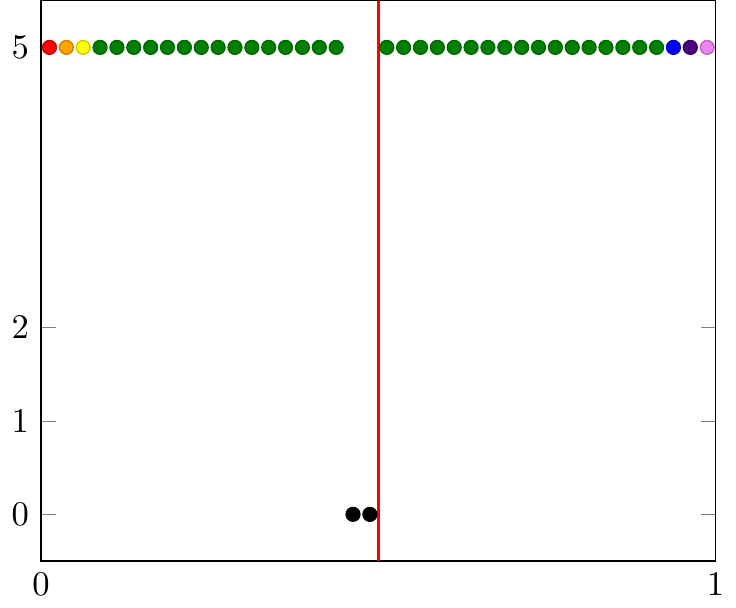} &
\includegraphics[width=0.3\textwidth,align=c]{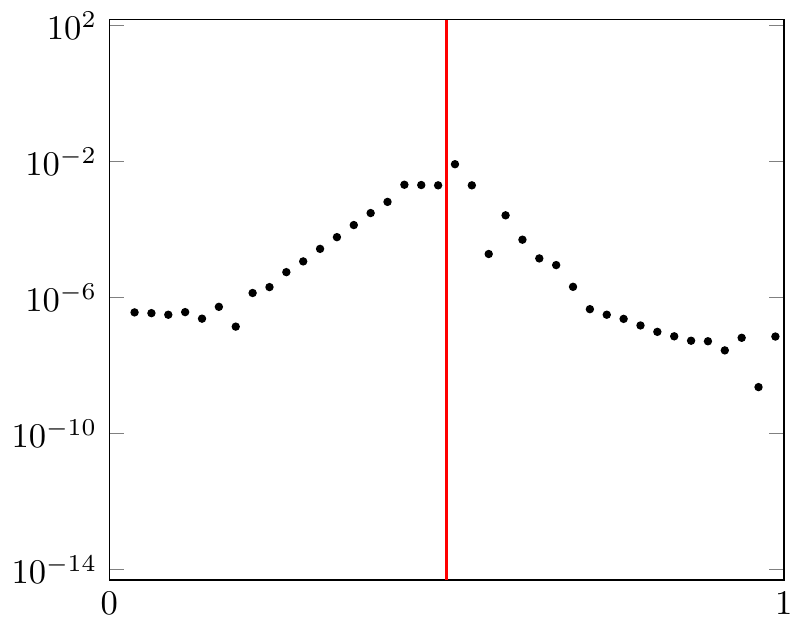}\\
AS & \includegraphics[width=0.3\textwidth,align=c]{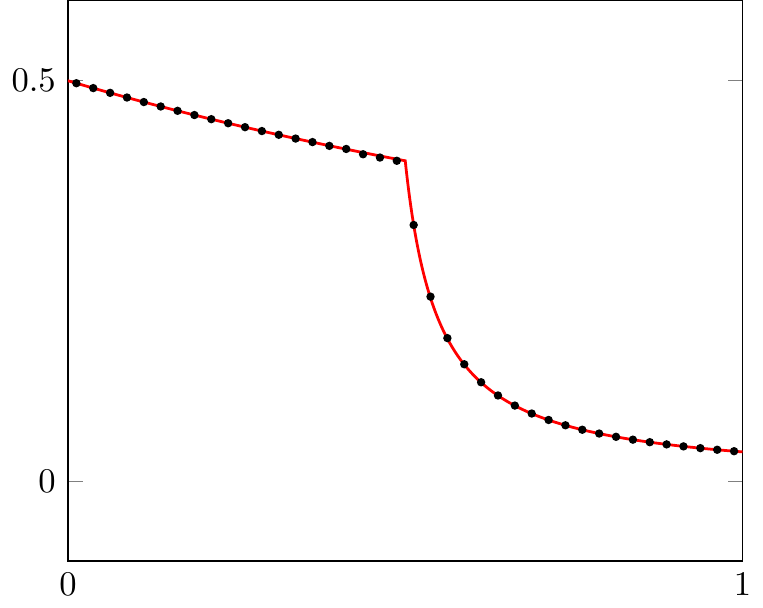} &
\includegraphics[width=0.3\textwidth,align=c]{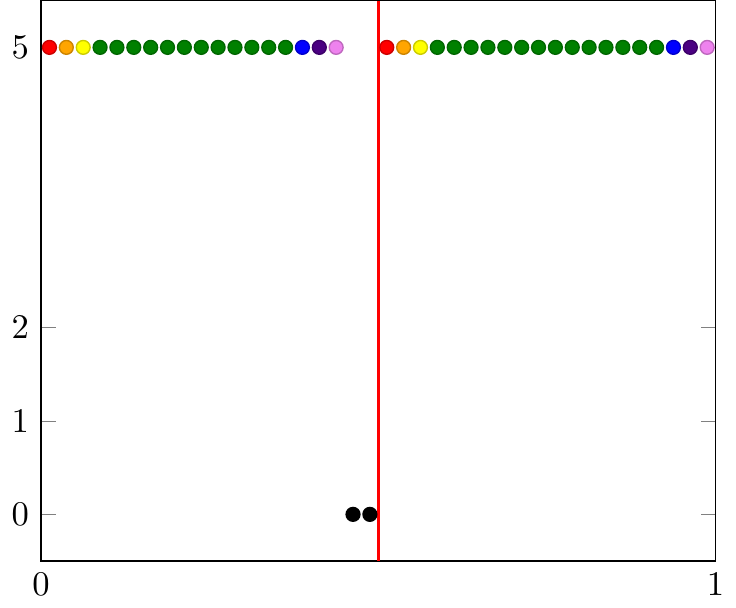} &
\includegraphics[width=0.3\textwidth,align=c]{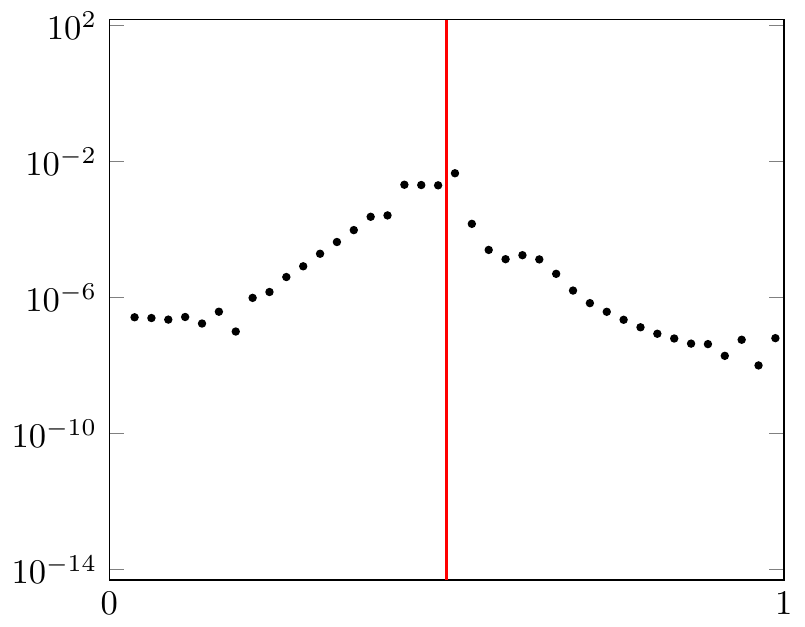}\\\bottomrule
\end{tabular}
\caption{Advection example 2 (irregular): TM1 solver for $I=40$.}
\label{fig:L_ex2_TM1}
\end{figure}

\begin{figure}\centering
\begin{tabular}{@{}cc@{}c@{}c@{}}\toprule
& exact and numerical solutions & CPD map & errors\\\midrule
MOOD & \includegraphics[width=0.3\textwidth,align=c]{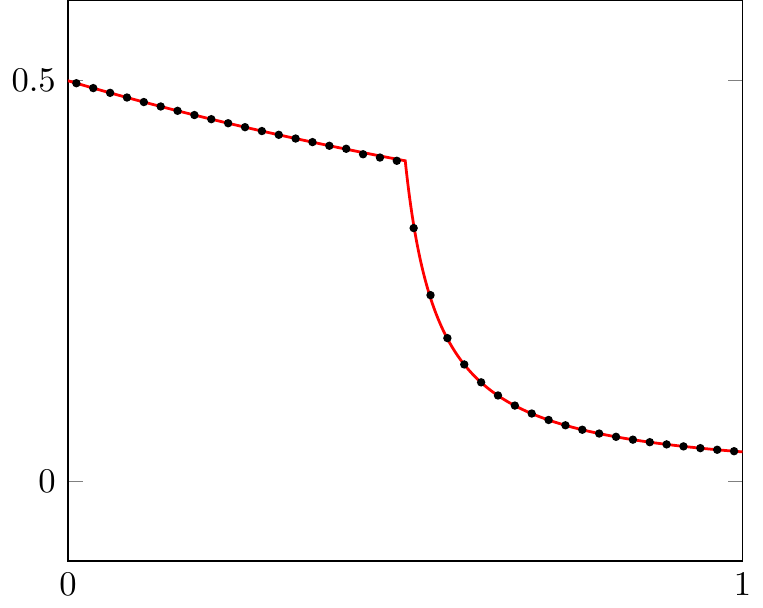} &
\includegraphics[width=0.3\textwidth,align=c]{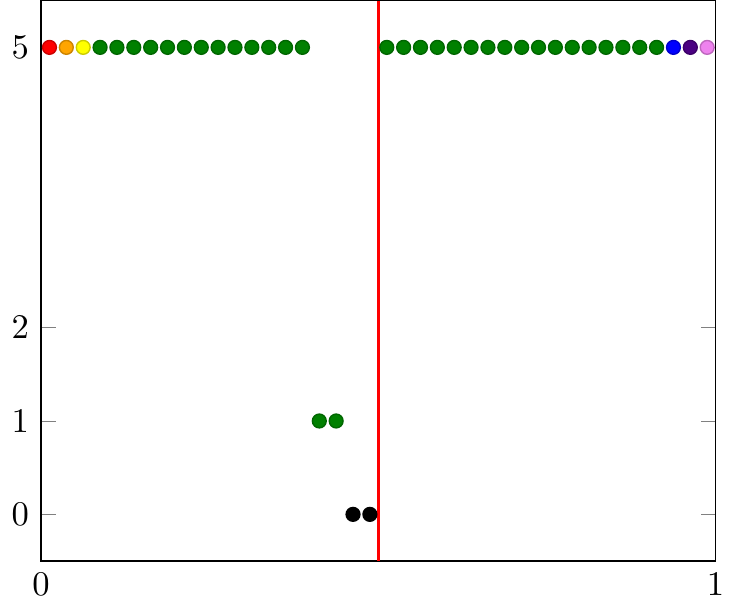} &
\includegraphics[width=0.3\textwidth,align=c]{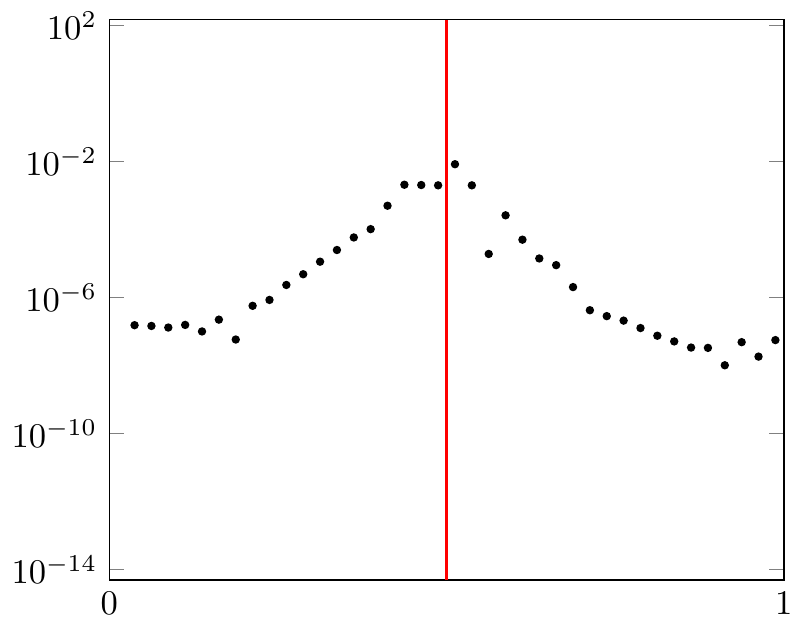}\\
AS & \includegraphics[width=0.3\textwidth,align=c]{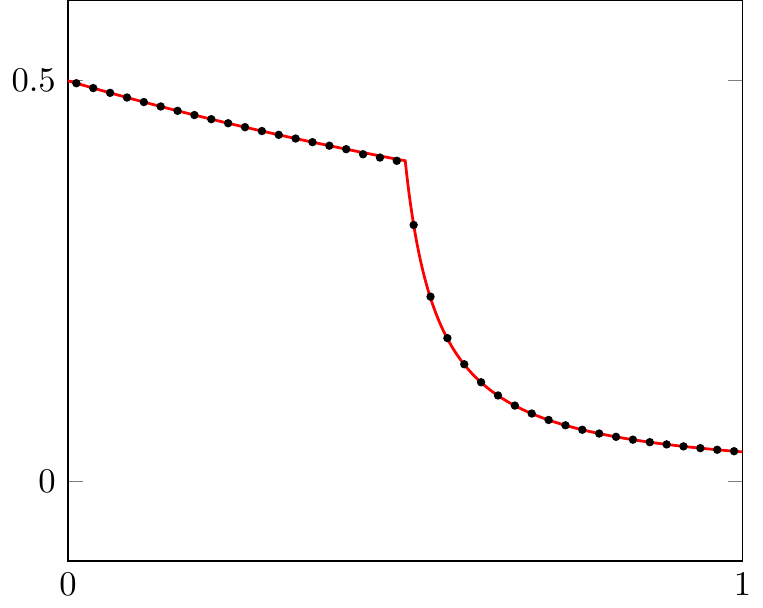} &
\includegraphics[width=0.3\textwidth,align=c]{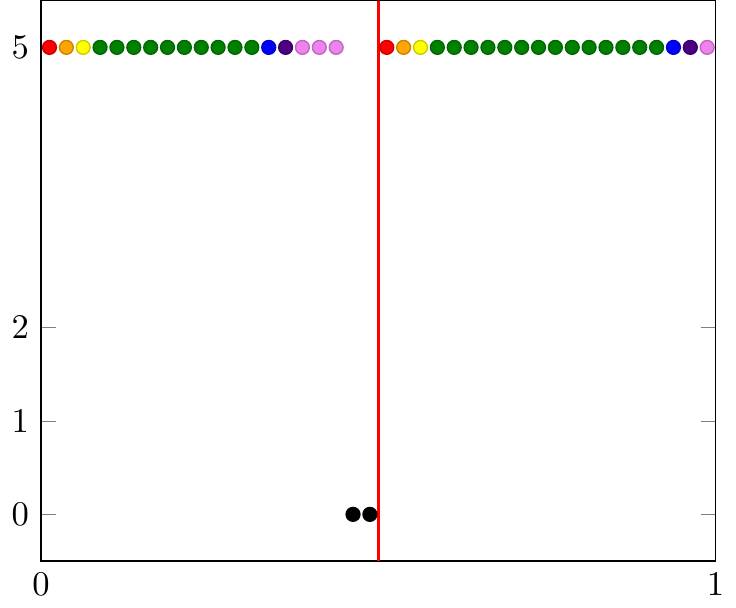} &
\includegraphics[width=0.3\textwidth,align=c]{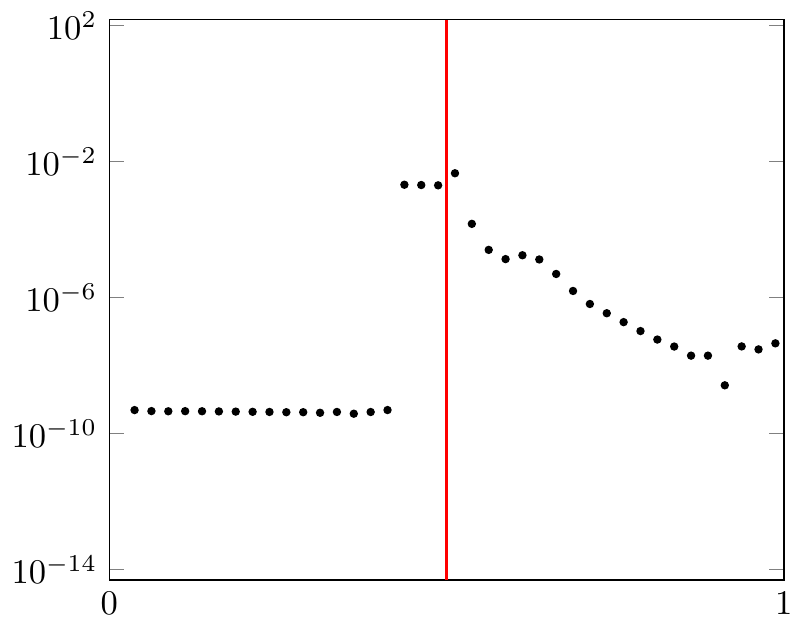}\\\bottomrule
\end{tabular}
\caption{Advection example 2 (irregular): TM2 solver for $I=40$.}
\label{fig:L_ex2_TM2}
\end{figure}

We have carried out the same simulation with the two Time Marching strategies 
%we reproduce 
in Figs.~\ref{fig:L_ex2_TM1} and \ref{fig:L_ex2_TM2} for TM1 and TM2, respectively. Curiously, the TM1 strategy produces a better CPD map in the case of a Static Stencil leading to a slightly better solution. 
The TM2 method exactly provides the same CPD map than the L method and, once again, the Adaptative Stencil strategy achieves a visible improvement of the polynomial reconstruction by selecting the best possible stencil. 

\subsection{B\"urgers' equation} \label{ssec:burgers}

Inviscid B\"urgers steady-state problem is an attractive benchmark since 
%we are dealing 
it deals with a simple scalar but non-linear model which may, nonetheless, generate discontinuous solutions. It reads
\[
\frac{\text{d}}{\text{d}x}\Big(\frac{\phi^2}{2}\Big)=S(x),\quad x\in \Omega,
\]
where $S\equiv S(x)$ is a given source term.
Dirichlet boundary conditions $\phi_{\text{L}}$ and $\phi_{\text{D}}$ are prescribed on the left and right boundaries of $\Omega$, respectively.

The numerical benchmark is the following one (see \cite{burgers_ex_irregular}):
\[
\begin{cases}
\frac{\partial\phi}{\partial t}+\frac{\partial}{\partial x}\Big(\frac{\phi^2}{2}\Big)=\sin(x)\cos(x), & x\in [0,\pi],\\
\phi(x,0)=\beta\sin(x), & x\in [0,\pi],\\
\phi(0,t)=\phi(\pi,t)=0, & t\in\mathds{R},
\end{cases}
\]
for which the steady-state solution is
\begin{itemize}
\item if $\beta\in]-\infty,-1]$: $\phi(x,\infty)=\sin(x)$;
\item if $\beta\in]-1,1[$: $\phi(x,\infty)=\begin{cases}\sin(x), & \text{if } x\in[0,x^\ast[\\-\sin(x), & \text{if } x\in[x^\ast,\pi],\end{cases}$
with $x^\ast=\pi-\arcsin\left(\sqrt{1-\beta^2}\right)$;
\item if $\beta\in[1,\infty[$: $\phi(x,\infty)=-\sin(x)$.
\end{itemize}

Note that the location of the shock is determined by parameter $\beta$ through the conservation of the initial total mass as
$$
\int_0^\pi\phi(x,t)\,\text{d}x=\int_0^\pi\beta\sin(x)\,\text{d}x=2\beta.
$$ 
%In our case we 
We take $x^\ast=\frac{3\pi}{4}$ which implies $\beta=\frac{\sqrt{2}}{2}$.
The numerical flux $\mathcal F(\phi_{i+\frac{1}{2},-},\phi_{i+\frac{1}{2},+};x_{i+\frac{1}{2}})$ is the one from Rusanov 
\[
\mathcal{F}(\phi_{i+\frac{1}{2},-},\phi_{i+\frac{1}{2},+};x_{i+\frac{1}{2}})=\frac{1}{2}\left(\frac{\phi_{i+\frac{1}{2},-}^2}{2}+\frac{\phi_{i+\frac{1}{2},+}^2}{2}\right)-\lambda_{i+\frac{1}{2}}(\phi_{i+\frac{1}{2},+}-\phi_{i+\frac{1}{2},-}),
\]
where $\lambda_{i+\frac{1}{2}}=\max(|\phi_i|,|\phi_{i+1}|)$.
Notice that the same detectors' chain used for linear problem is also employs in the B\"urger's case. 

\begin{figure}\centering
\begin{tabular}{@{}cc@{}c@{}c@{}}\toprule
& exact and numerical solutions & CPD map & errors\\\midrule
MOOD & \includegraphics[width=0.3\textwidth,align=c]{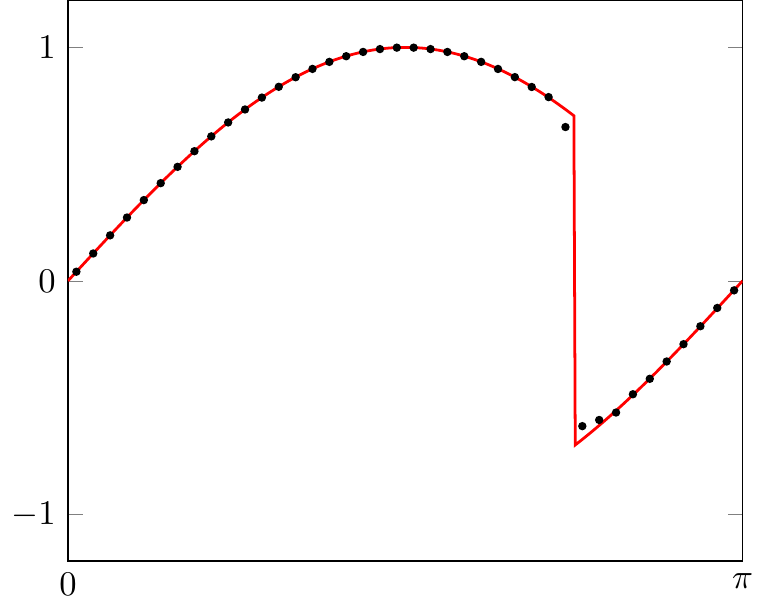} &
\includegraphics[width=0.3\textwidth,align=c]{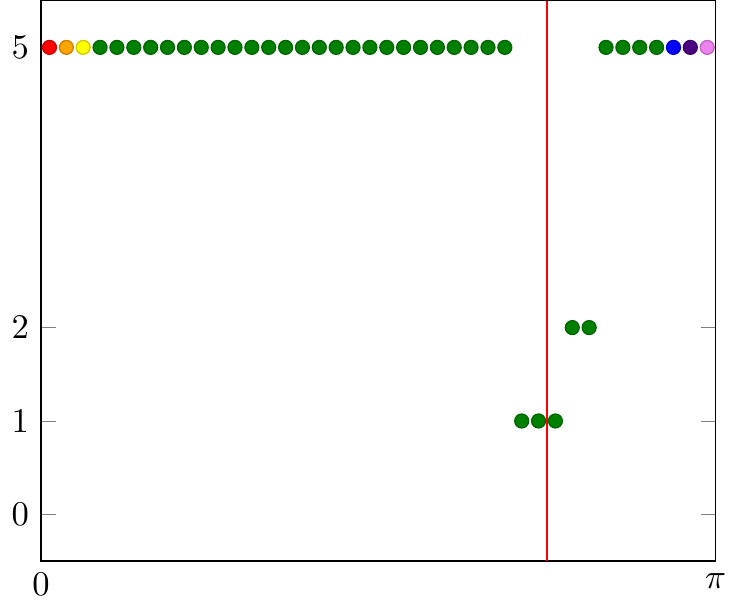} &
\includegraphics[width=0.3\textwidth,align=c]{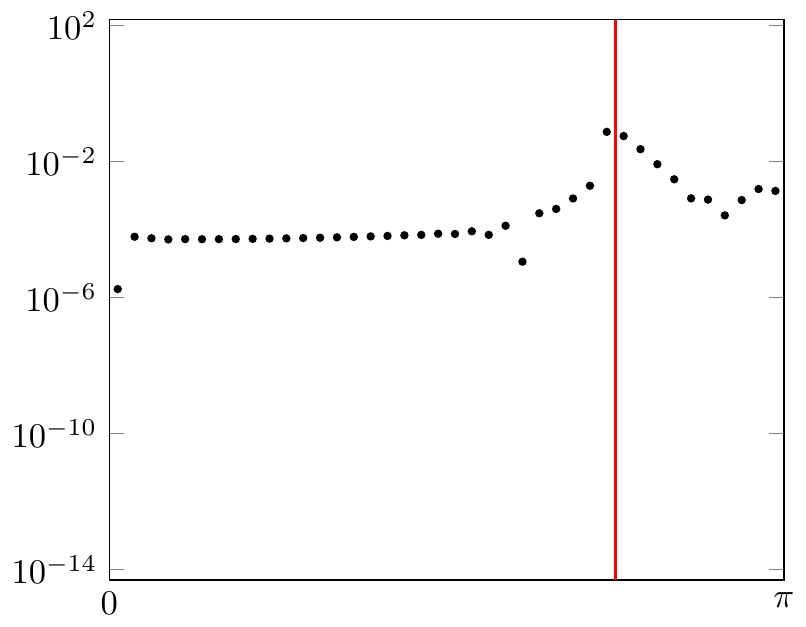}\\
AS & \includegraphics[width=0.3\textwidth,align=c]{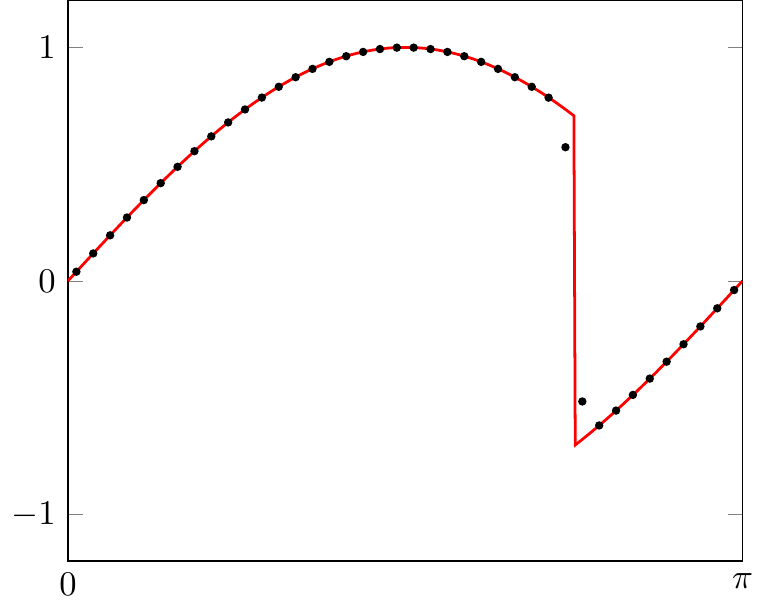} &
\includegraphics[width=0.3\textwidth,align=c]{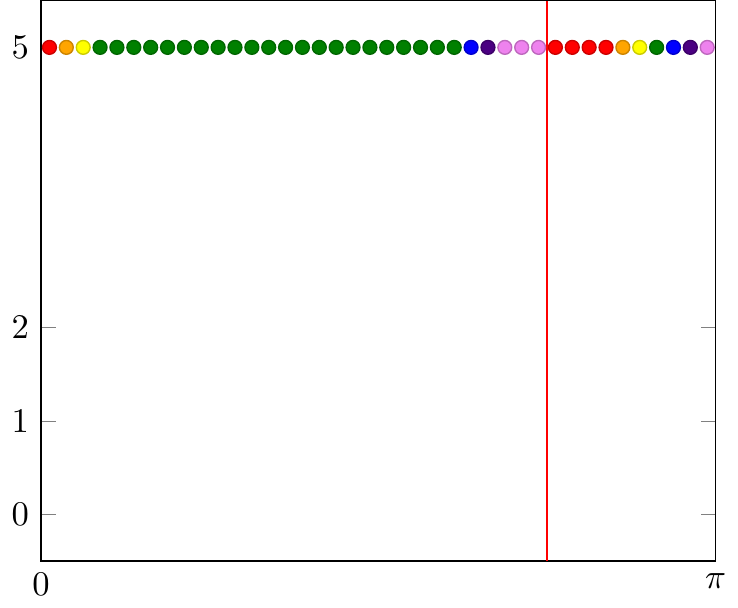} &
\includegraphics[width=0.3\textwidth,align=c]{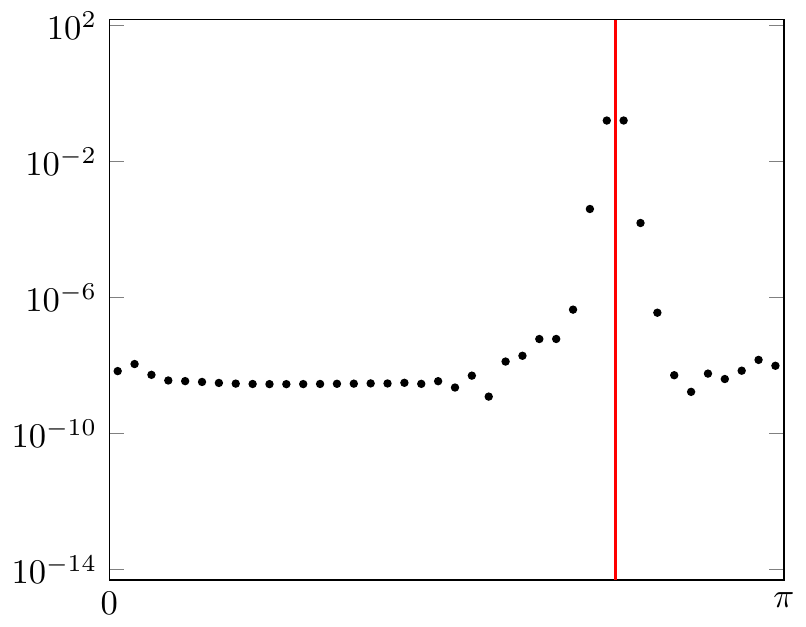}\\\bottomrule
\end{tabular}
\caption{\label{fig:burger_benchmark} B\"urgers example: NL solver for $I=40$.}
\end{figure}

The first simulation has been carried out with the original MOOD method using the Static Stencil (first row in Fig.~\ref{fig:burger_benchmark}) whereas we present, in the second row, the Adaptative Stencil re\-cons\-truc\-tion to improve the accuracy by preserving the high degree of the polynomial reconstruction. We note that no oscillations have been created in both cases. The CPD map with the colors' code indicate that the Adaptative Stencil strategy manages to recover the optimal degree by appropriately upwinding the stencils. The last column confirms the spectacular gain because the error is lowered by several orders of magnitude in the shock vicinity. 
%The chief point 
The direct consequence is an excellent improvement of the solution quality in regular zones.

\begin{figure}\centering
\begin{tabular}{@{}cc@{}c@{}c@{}}\toprule
& exact and numerical solutions & CPD map & errors\\\midrule
MOOD & \includegraphics[width=0.3\textwidth,align=c]{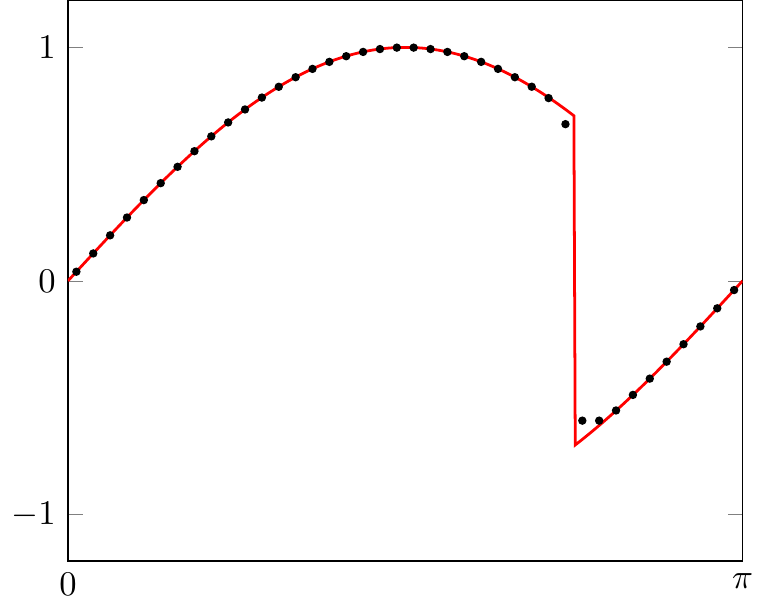} &
\includegraphics[width=0.3\textwidth,align=c]{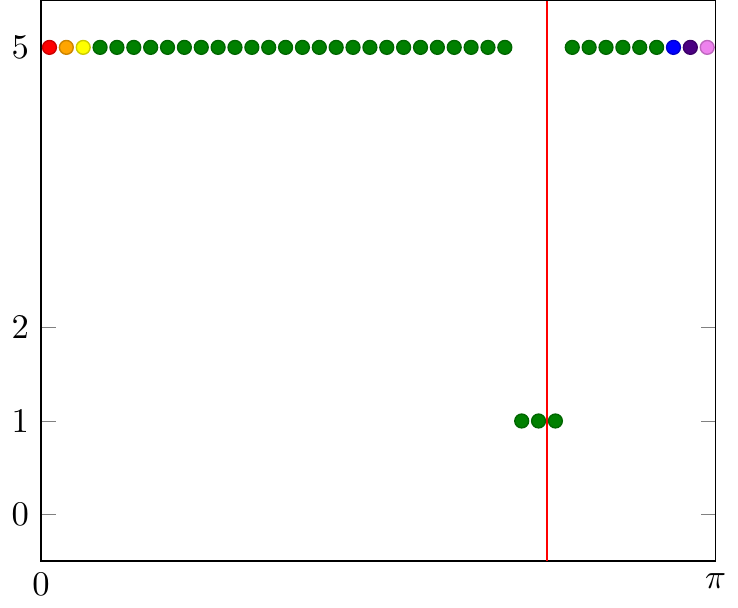} &
\includegraphics[width=0.3\textwidth,align=c]{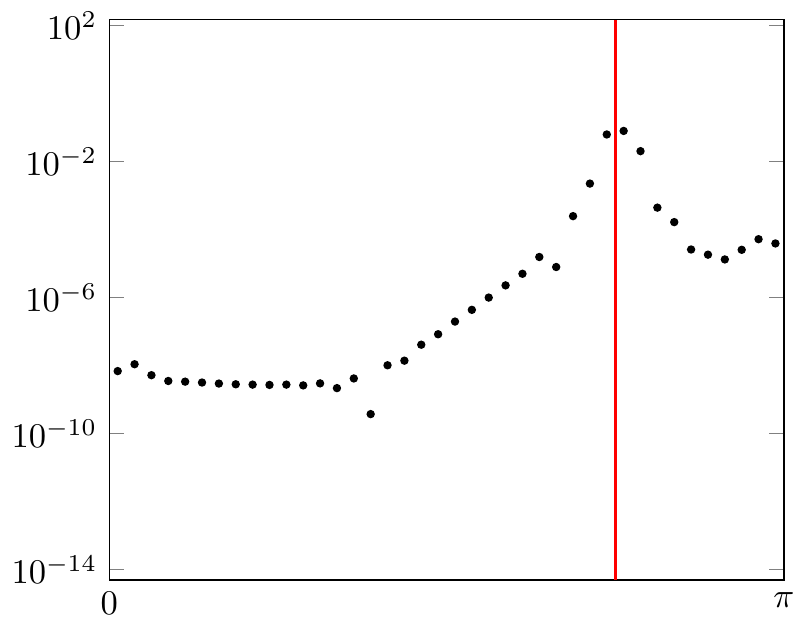}\\
AS & \includegraphics[width=0.3\textwidth,align=c]{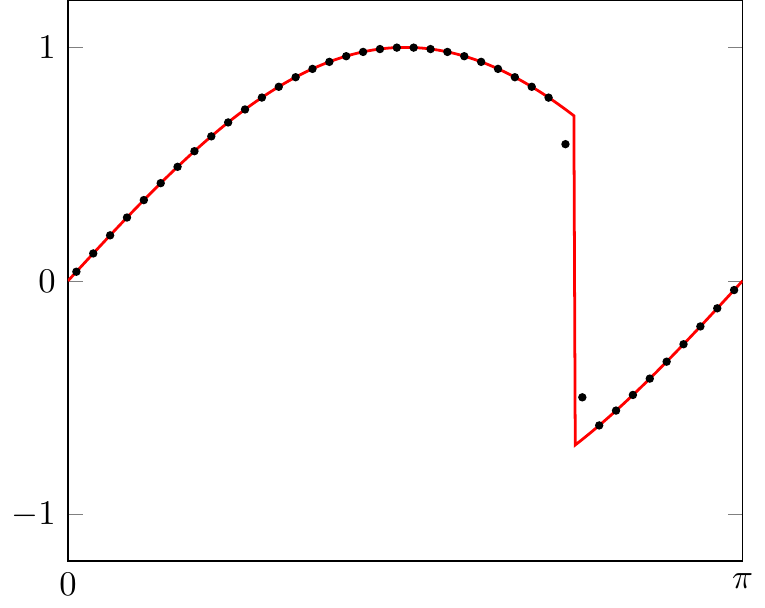} &
\includegraphics[width=0.3\textwidth,align=c]{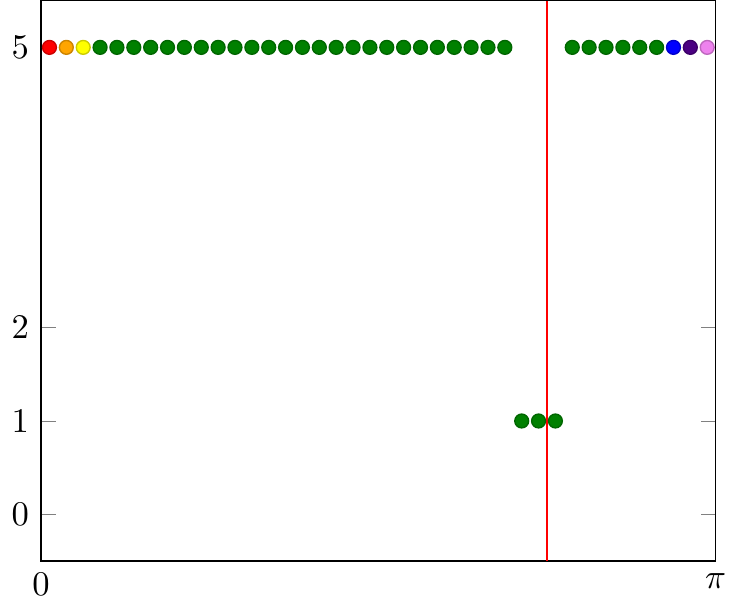} &
\includegraphics[width=0.3\textwidth,align=c]{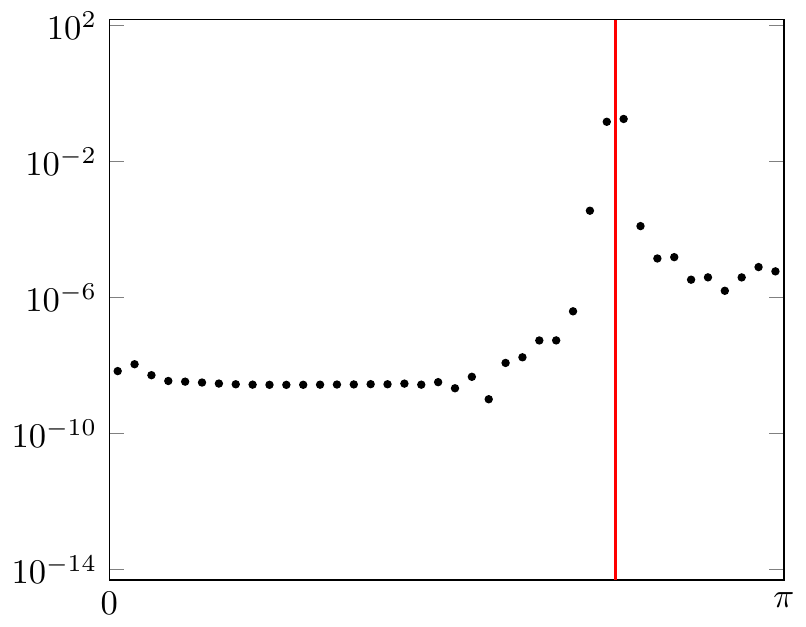}\\\bottomrule
\end{tabular}
\caption{\label{fig:burger_benchmark_TM1} B\"urgers example: TM1 solver for $I=40$.}
\end{figure}

\begin{figure}\centering
\begin{tabular}{@{}cc@{}c@{}c@{}}\toprule
& exact and numerical solutions & CPD map & errors\\\midrule
MOOD & \includegraphics[width=0.3\textwidth,align=c]{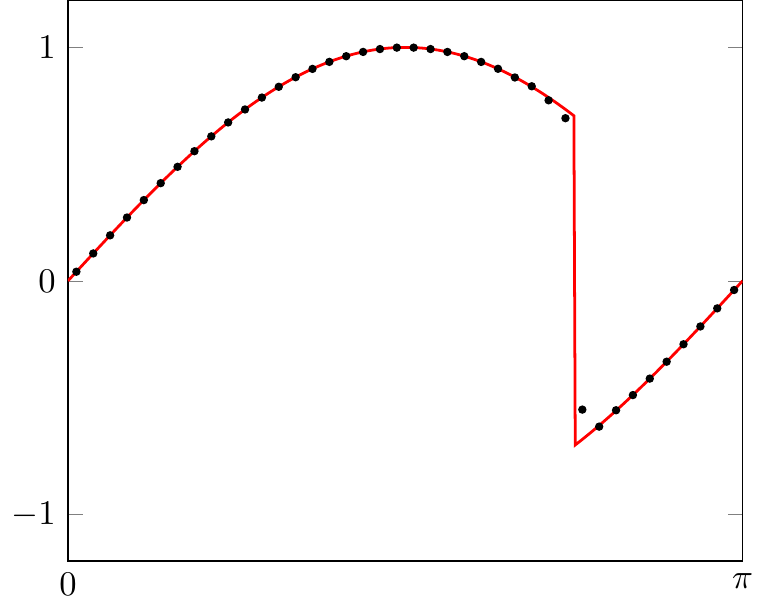} &
\includegraphics[width=0.3\textwidth,align=c]{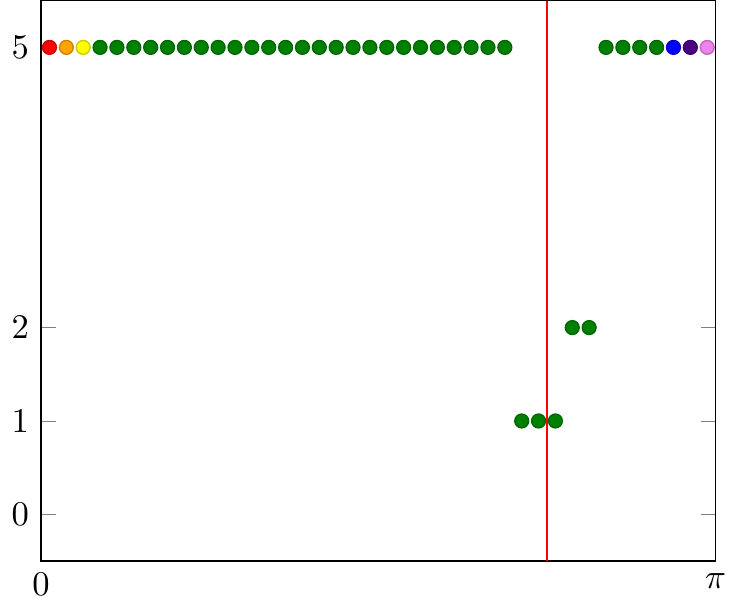} &
\includegraphics[width=0.3\textwidth,align=c]{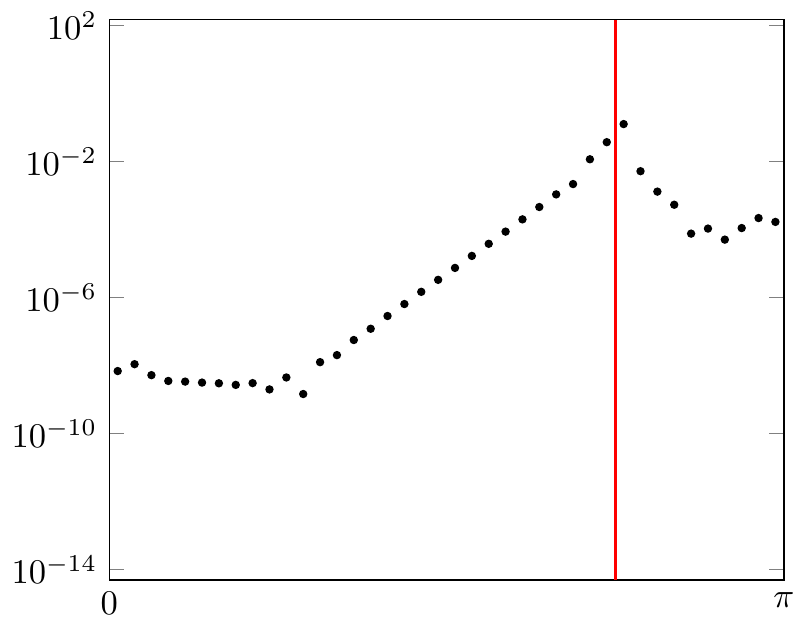}\\
AS & \includegraphics[width=0.3\textwidth,align=c]{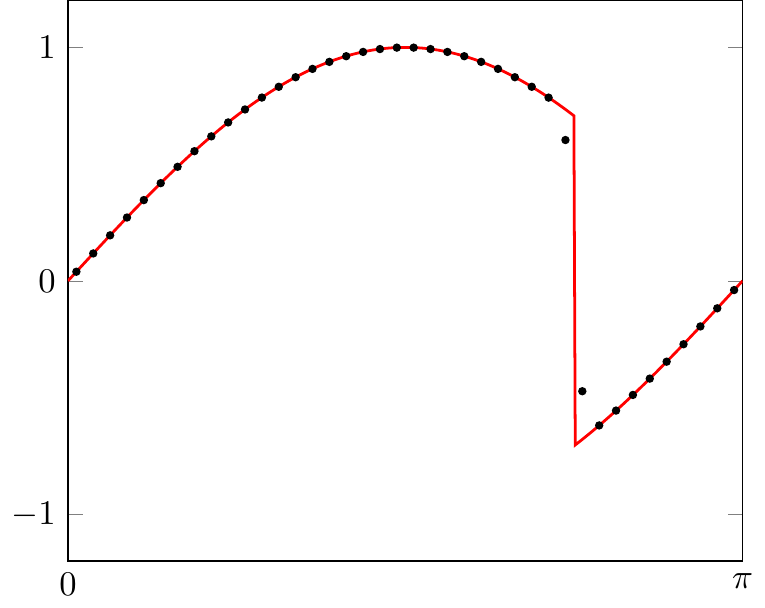} &
\includegraphics[width=0.3\textwidth,align=c]{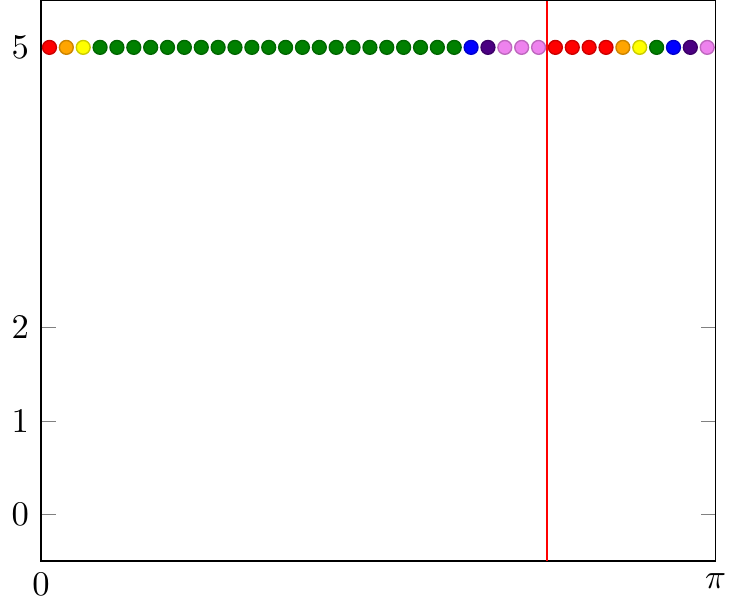} &
\includegraphics[width=0.3\textwidth,align=c]{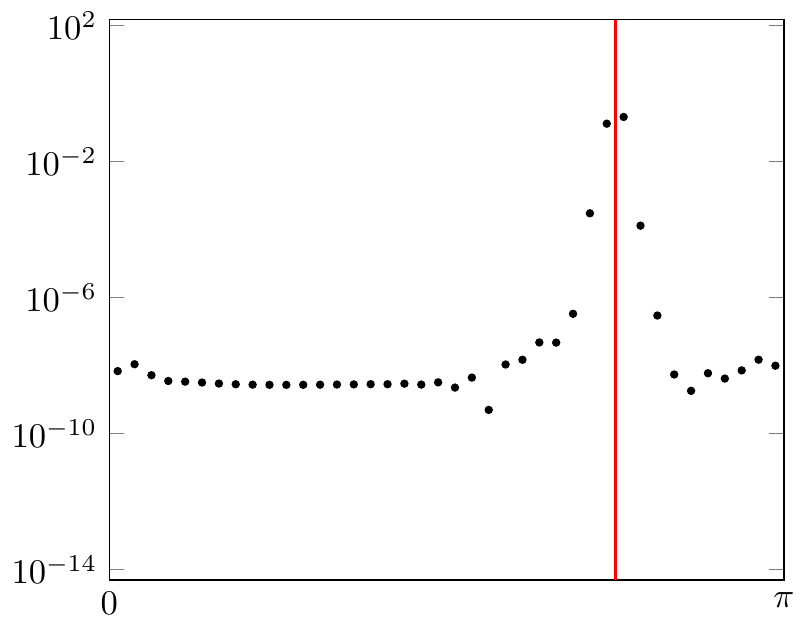}\\\bottomrule
\end{tabular}
\caption{\label{fig:burger_benchmark_TM2} B\"urgers example: TM2 solver for $I=40$.}
\end{figure}

We can draw similar conclusions with the TM1 and TM2 method 
%and display in
in view of Figs.~\ref{fig:burger_benchmark_TM1} and \ref{fig:burger_benchmark_TM2}. %the approximated solution, the CPD map together with the colors'' code and the errors.

\subsection{Euler system of equations} \label{ssec:euler}

Let $\bm{U}=(\rho,\rho u ,E)^T$ be the mass, momentum, and total energy conserved quantities and $\bm{V}=(\rho,u,p)^T$, the vector of primitive/state variables: density, velocity, and pressure. We seek a solution of the steady-state 1D Euler equations with source terms given by
\begin{equation}
\label{eq:Euler1D}
\frac{\textrm{d}\bm{F}(\bm{U})}{\textrm{d}x}=\bm{S}, \qquad\text{in }\Omega,
\end{equation}
where the flux vector is given by $\bm{F}(\bm{U})=(\rho u,\rho u^2+p,u(E+p))^T$ and the regular vector source term is given by $\bm{S}\equiv\bm{S}(x)=(\frac{\textrm{d}D}{\textrm{d}x},\frac{\textrm{d}F}{\textrm{d}x},\frac{\textrm{d}H}{\textrm{d}x})^T\equiv(\frac{\textrm{d}D}{\textrm{d}x}(x),\frac{\textrm{d}F}{\textrm{d}x}(x),\frac{\textrm{d}H}{\textrm{d}x}(x))^T$.
Total energy $E$, specific internal energy $e$, and sound-speed $a$ are 
%expressed as
related to one another by
\[
E=\frac{1}{2}\rho u^2+\rho e, \quad
e=\frac{p}{\rho(\gamma-1)}, \quad
a=\sqrt{\frac{\gamma p}{\rho}},
\]
respectively, where $\gamma>1$ is the ratio of specific heat. System~\eqref{eq:Euler1D} is completed with Dirichlet boundary conditions
\[
\bm{V}(x_{\text{L}})=\bm{V}_\text{L}\equiv(\rho_\text{L},u_\text{L},p_\text{L})^T
\quad\text{ and }\quad
\bm{V}(x_{\text{R}})=\bm{V}_\text{R}\equiv(\rho_\text{R},u_\text{R},p_\text{R})^T.
\]
We deduce the algebraic system to be solved 
\begin{equation} 
\label{eq:euler_cond}
%	\label{eq:D}
	\rho u = D,    \;
%	\label{eq:F}
	\rho u^2+p = F, \;
%	\label{eq:H}
	u(E+p) = H,\;
%	\label{eq:E} 
	\frac{1}{2}\rho u^2+\frac{p}{\gamma-1}=E.
\end{equation}
Notice that $D$, $F$, and $H$ depend on additional constants we shall fix with the boundary conditions.
Analytic solutions for the steady-state problem are detailed in \cite{1D_SS}. 

Let us denote by $\bm{\mathcal U}=(\bm{U}_i)_{i=1,\ldots,I}$ the
$3\times I$ matrix containing the approximations of the mean-values for $\rho$, $\rho u$, and $E$ while
$\bm{\mathcal U}^k$ stands for the data associated to stage $k$. The numerical flux 
$\mathcal F (\bm{U}_-,\bm{U}_+;x)$ is chosen to be the HLL/HLLE one proposed in \cite{hll} where $\bm{U}_-$ and $\bm{U}_+$ are the left and right states at interface $x$
\[
\mathcal{F}(\bm{U}_-,\bm{U}_+;x)=
\begin{cases} 
\bm{U}_- & \text{if $s_-\geq 0$},\\
\frac{s_+\bm{F}(\bm{U}_-)-s_-\bm{F}(\bm{U}_+)+s_+s_-(\bm{U}_+-\bm{U}_-)}{s_+-s_-} & \text{if $s_-\leq 0 \leq s_+$} \, ,\\
\bm{U}_+ & \text{if $s_+\leq 0$},
\end{cases}
\]
where we %take the numerical velocities 
consider the wave speeds to be $s_-=\min(u_--a_-,u_+-a_+)$ and $s_+=\max(u_-+a_-,u_++a_+)$ for the sake of simplicity. 
\begin{remI}
	Notice that such a flux does not strictly preserve the positivity of the density nor pressure. If one requires a formally positivity preserving flux, the Rusanov one equipped with the adequate numerical velocity should be preferred (see for instance \cite{munz91} for details on positivity preservation).
\end{remI}

For given CPD and CS maps $\mathcal M^k,\mathcal S^k$, we introduce the nonlinear operator
\[
\bm{\mathcal U}\in\mathbb R^{3\times I}\to \mathcal{G}(\bm{\mathcal U},\mathcal M^k,\mathcal S^k) \in\mathbb R^{3\times I}
\]
and we seek for $\bm{\mathcal U}^{k+1}$ such that  $\mathcal{G}(\bm{\mathcal U}^{k+1},\mathcal M^k,\mathcal S^k)=0$.

\begin{figure}
	\centering
	\hspace{-0.75cm} \includegraphics[width=1.0\textwidth]{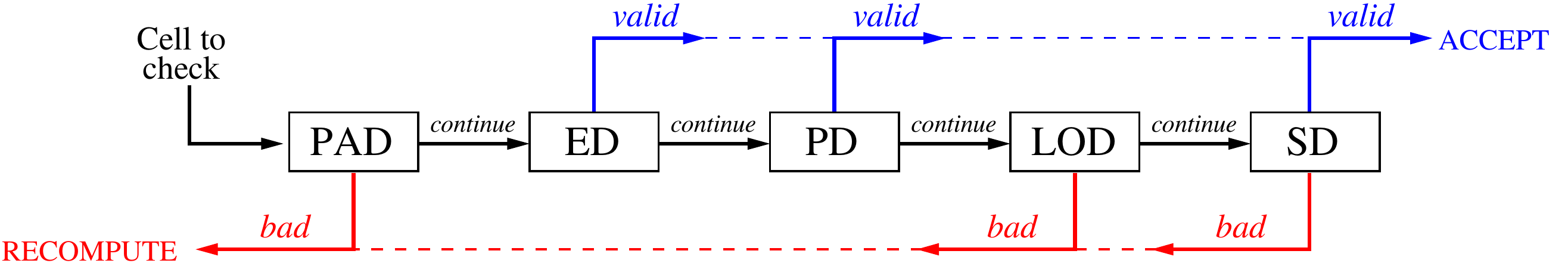}
	\caption{Chain of detectors used for Euler system of equations
		to check if a cell value is valid and accepted, or, if it is discarded 
		and need to be recomputed after reduction of the local cell polynomial degree.
		PAD: Physical Admissible Detector,
		ED: Extrema Detector,
		PD: Plateau Detector,
		LOD: Local Oscillation Detector,
		SD: Smoothness Detector.
	}
	\label{fig:chain_detector_euler}
\end{figure}

%Again we consider either a nonlinear solver (NL) and the time-marching (TM) technique to solve this set of equations. 
%Notice that for the NL case, an additional condition is considered to fix the shock 
%position when dealing with two different regimes on the left and right boundaries.
%For the TM method, the initial condition for density imposes the position of the shock.
% detector chain and cascade
The detector chain associated to the Euler equation follows the chain depicted in Fig.~\ref{fig:chain_detector_euler}. Here the PAD box checks the positivity of density and pressure.
% --- FIG --- detector chain
% --- FIG ---
Because we consider only one polynomial degree per cell for each component of $\bm{\mathcal U}$, \textit{i.e.}, density, momentum, and total energy, we determine only one polynomial degree map. Only the density variable is tested for the ED, LOD, and SD detection criteria and the CS map is determined accordingly.
The MOOD loop is escaped if $\mathcal{M}^{k+1}=\mathcal{M}^{k}$.

\begin{figure}
	\centering
	\includegraphics{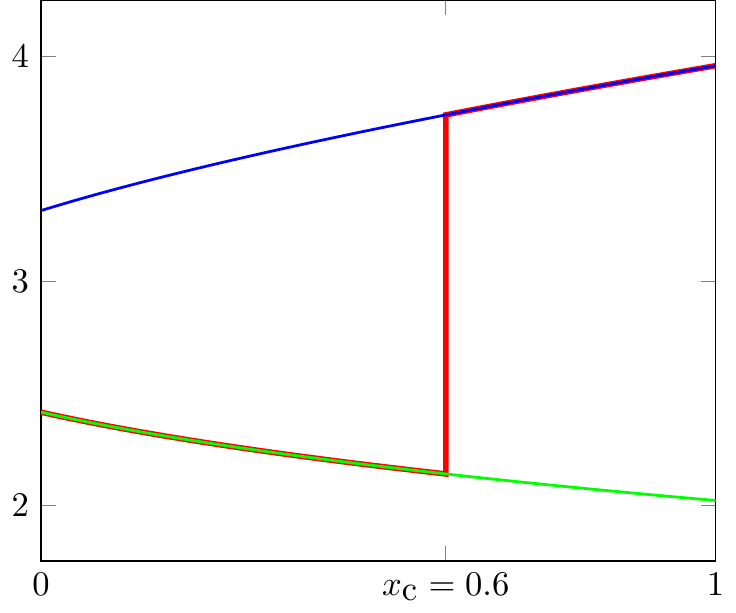}	
	\caption[]{Plots of the supersonic and subsonic branches $\rho_\textup{sup}$ (\rhosup) and $\rho_\textup{sub}$ (\rhosub) from~\eqref{eq:supsub},
		and the exact solution of the irregular density (shock located at $x_{\text{c}}=0.6$) employed in this paper (\exactsol).}
	\label{fig:EulerrhoSupAndrhoSub}
\end{figure}

To assess the validity of the method, we consider the following data. The domain is given by $x_{\text{L}}=0$, $x_{\text{R}}=1$ and we take $D(x)=1$, $F(x)=0.027x + 0.6137$, and $H(x)=0.375$.  The first two conditions in~\eqref{eq:euler_cond} for the admissible space are trivially satisfied. The third condition, $\frac{DH}{F^2}\leq\frac{\gamma^2}{2(\gamma^2-1)}$ is also, see \cite{1D_SS}.
%\[
%\max_{x\in[0,1]}\frac{DH}{F^2}=\max_{x\in[0,1]}\frac{0.375}{(0.027x + 0.6137)^2}=0.9957\leq\frac{\gamma^2}{2(\gamma^2-1)}=1.0208.
%\]
The expressions for the supersonic and the subsonic solution branches are given by
\begin{gather} \label{eq:supsub}
	\rho_\textup{sup}(x)=0.126x - 3.3333\sqrt{1.96(0.027x + 0.6137)^2 - 0.72}+ 2.8639,\\
	\rho_\textup{sub}(x)=0.126x + 3.3333\sqrt{1.96(0.027x + 0.6137)^2 - 0.72}+ 2.8639,
\end{gather}
respectively, and depicted in Fig.~\ref{fig:EulerrhoSupAndrhoSub}.
% --- FIG --- exact Euler figure
%   To create the fgure run gnuplot script in DATA/Euler/euler_ex.gnu

Prescribing a supersonic state on the left boundary and a subsonic one on the right produces a solution with a shock wave located at $x_\textup{c}$. 
%We now present this last situation.
In order to enforce the shock position at $x_\textup{c}=0.6$, we set  
\begin{equation} \label{eq:constraintA}
A=
\int_{x_{\textup{L}}}^{x_{\textup{R}}} \rho^0(x)\,\textup{d}x=
%\int_{x_{\textup{L}}}^{x_{\textup{R}}} \left((\rho_\text{L}-\rho_\text{R})(1-x)^{0.4550}+\rho_\text{R}\right)\,\textup{d}x =
2.8975. 
\end{equation}
for an initial density profile 
\begin{equation}
\rho^0(x) = 
\begin{cases}
\rho_\text{L} & \text{ if $x\leq x_0=0.6873$} \\
\rho_\text{R} & \text{ otherwise} \\
\end{cases}
\end{equation}
where $x_0$ is determined by solving~\eqref{eq:constraintA}. We initialize the other variables as $u^0=D/\rho^0$, $p^0$, and $E^0$ derive from relation~ \eqref{eq:euler_cond}. The exact solution for the density variable is depicted in red in Fig.~\ref{fig:EulerrhoSupAndrhoSub}.

In \cite{1D_SS} a numerical studied as been carried out to analyze the MOOD scheme behavior and we refer the reader to these results which can be summarized as follows: first, the unlimited high order schemes produce oscillatory results, and, second, the MOOD scheme is able to mitigate those spurious oscillations by reducing the CPD map in the vicinity of the shock position, \textit{i.e.}, we get low order polynomial reconstructions for the cells close to the discontinuity. 

\begin{figure}\centering
\begin{tabular}{@{}cc@{}c@{}c@{}}\toprule
& exact and numerical solutions & CPD map & errors\\\midrule
MOOD & \includegraphics[width=0.3\textwidth,align=c]{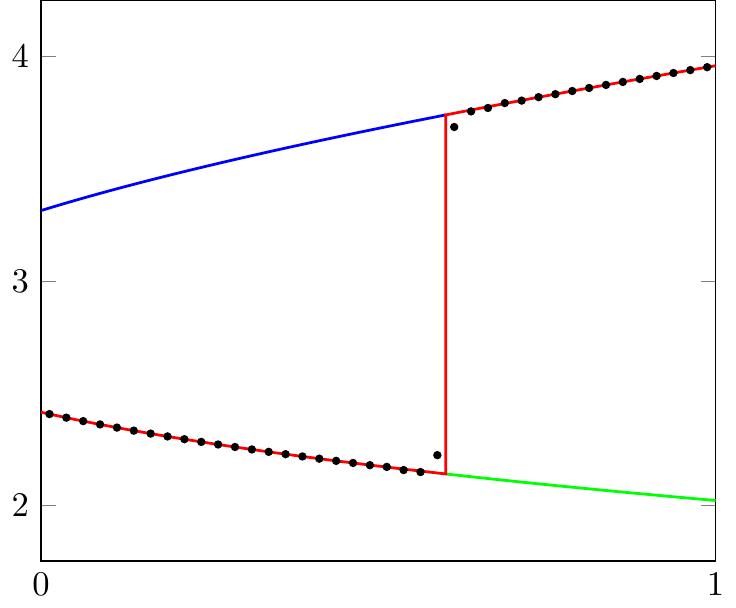} &
\includegraphics[width=0.3\textwidth,align=c]{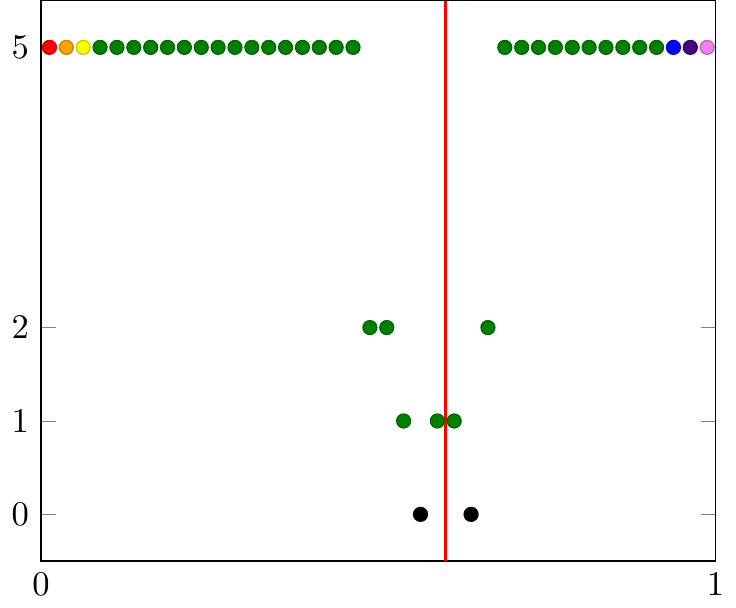} &
\includegraphics[width=0.3\textwidth,align=c]{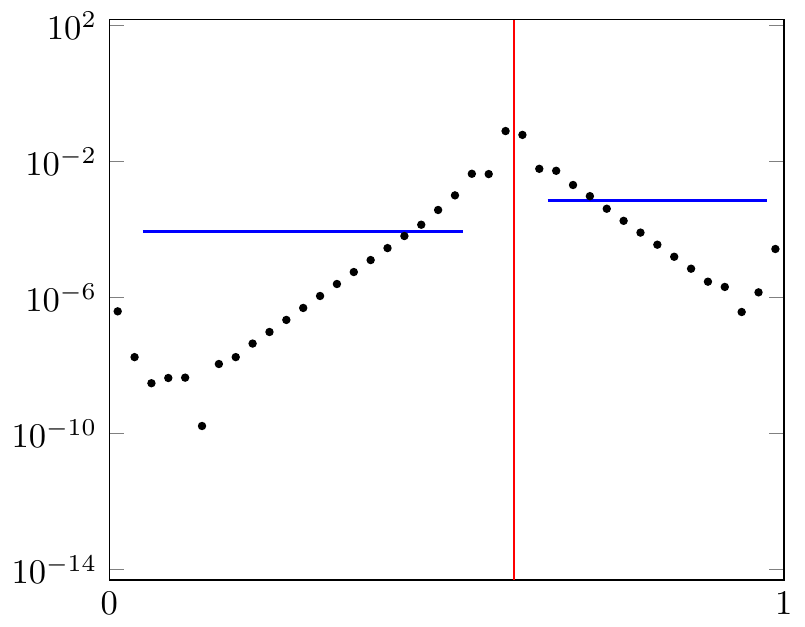}\\
AS & \includegraphics[width=0.3\textwidth,align=c]{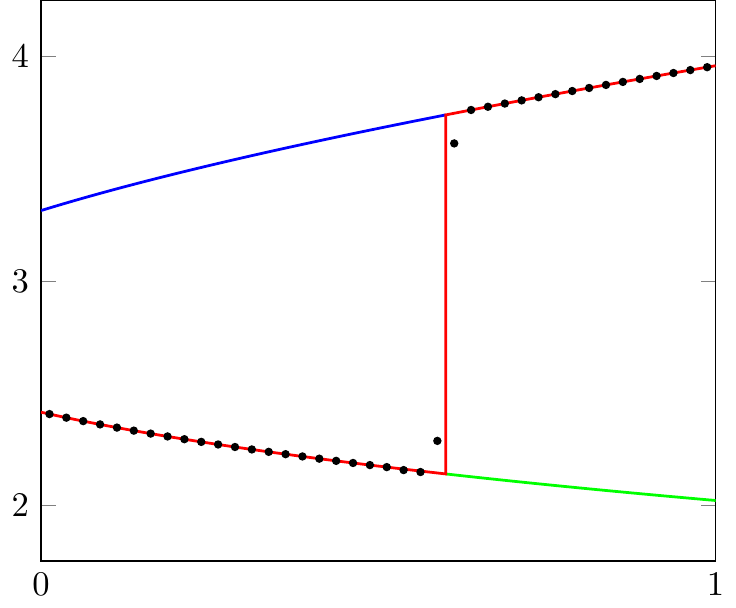} &
\includegraphics[width=0.3\textwidth,align=c]{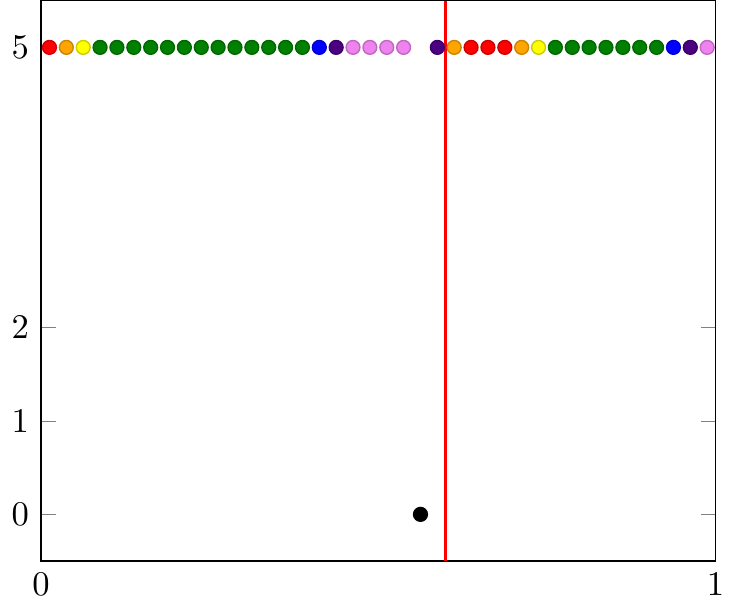} &
\includegraphics[width=0.3\textwidth,align=c]{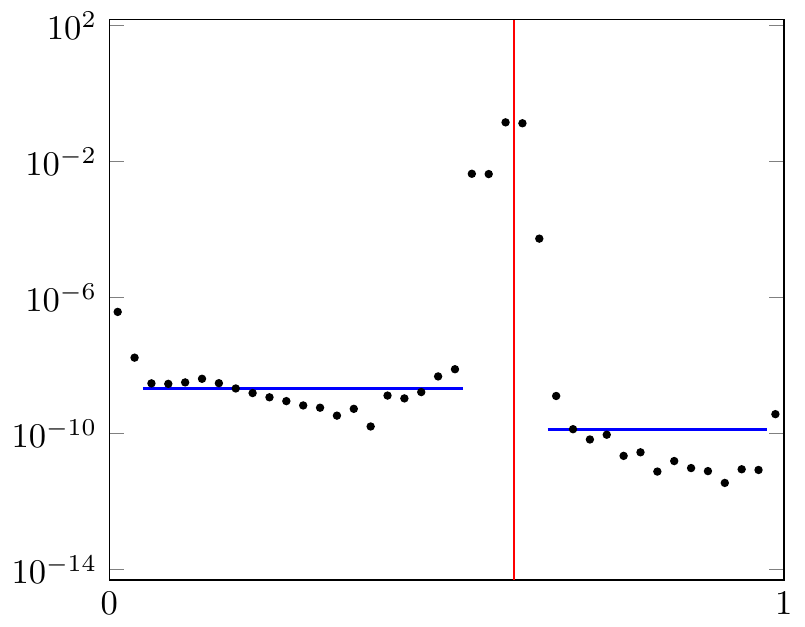}\\\bottomrule
\end{tabular}
\caption{\label{fig:euler_benchmarkI40} Euler example: NL solver for $I=40$ (shock location at an interface).}
\end{figure}

\begin{figure}\centering
\begin{tabular}{@{}cc@{}c@{}c@{}}\toprule
& exact and numerical solutions & CPD map & errors\\\midrule
MOOD & \includegraphics[width=0.3\textwidth,align=c]{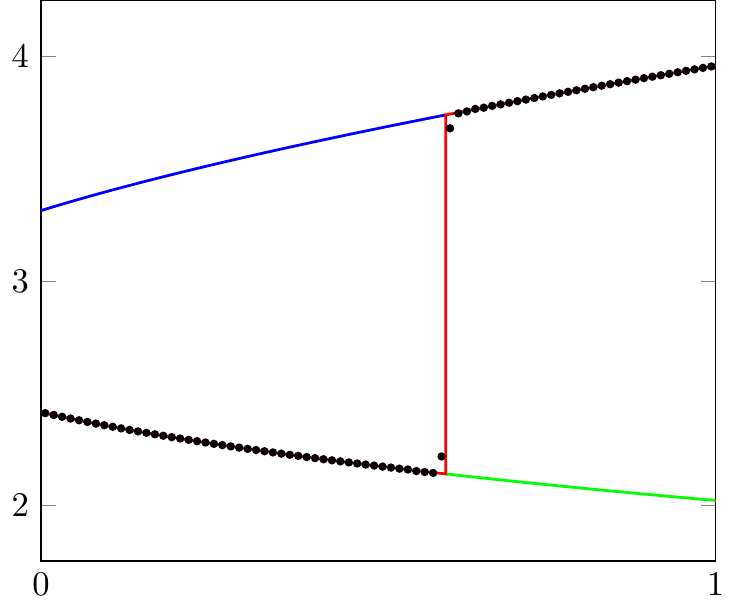} &
\includegraphics[width=0.3\textwidth,align=c]{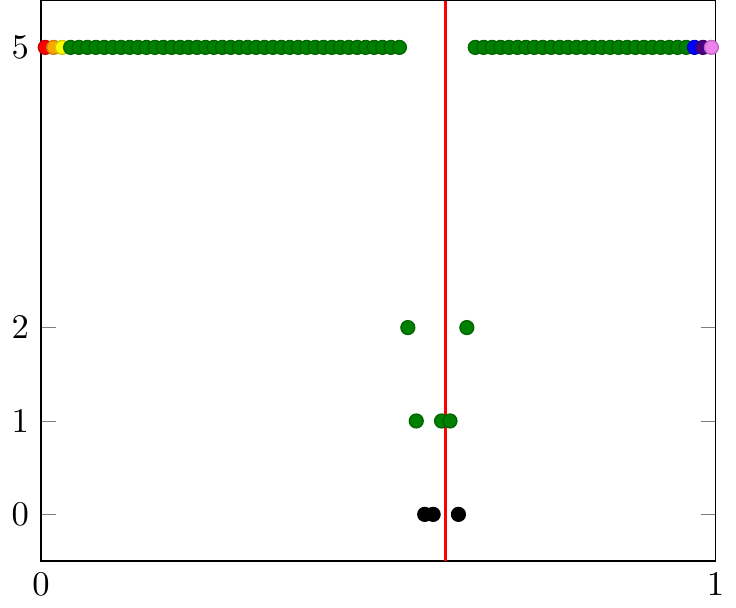} &
\includegraphics[width=0.3\textwidth,align=c]{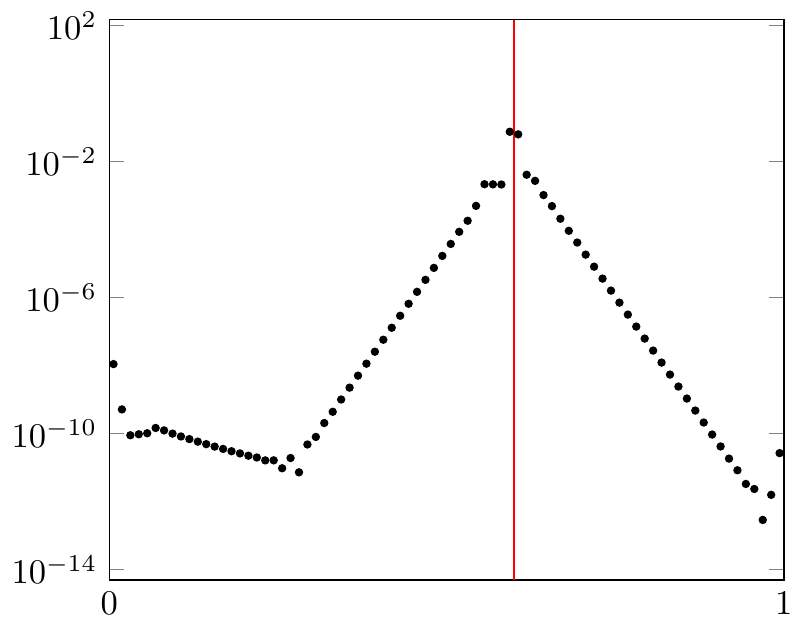}\\
AS & \includegraphics[width=0.3\textwidth,align=c]{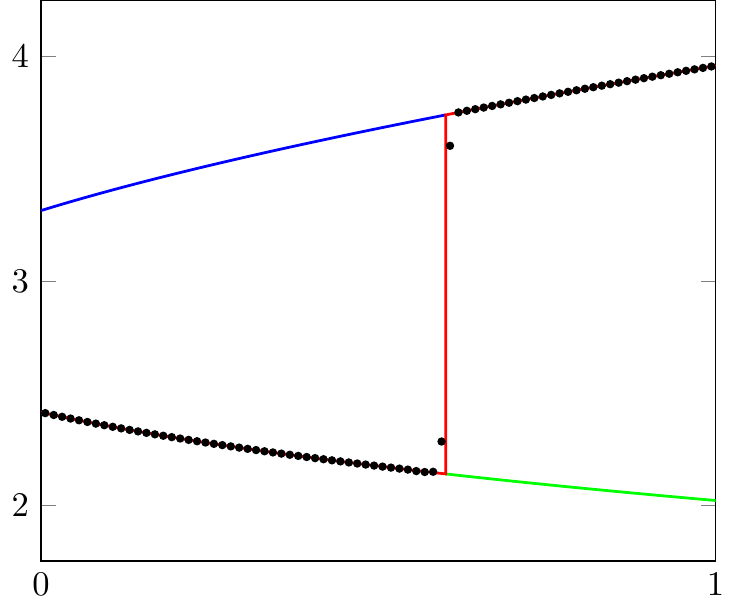} &
\includegraphics[width=0.3\textwidth,align=c]{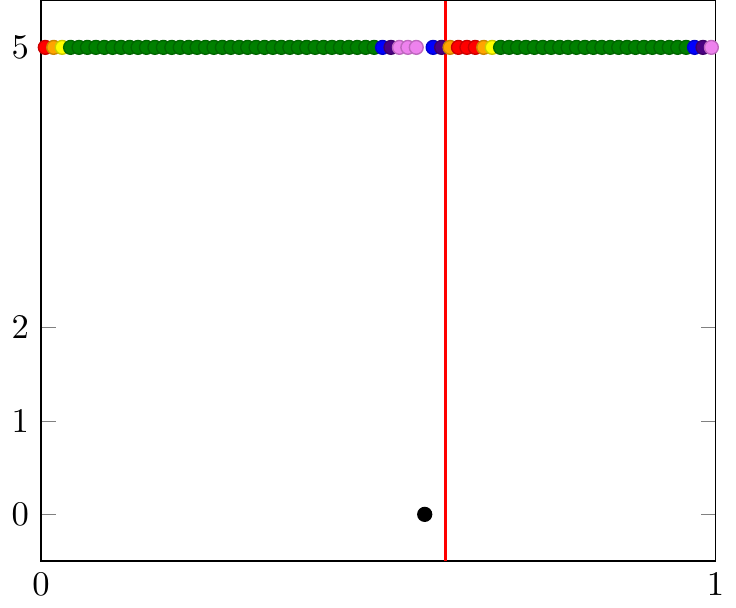} &
\includegraphics[width=0.3\textwidth,align=c]{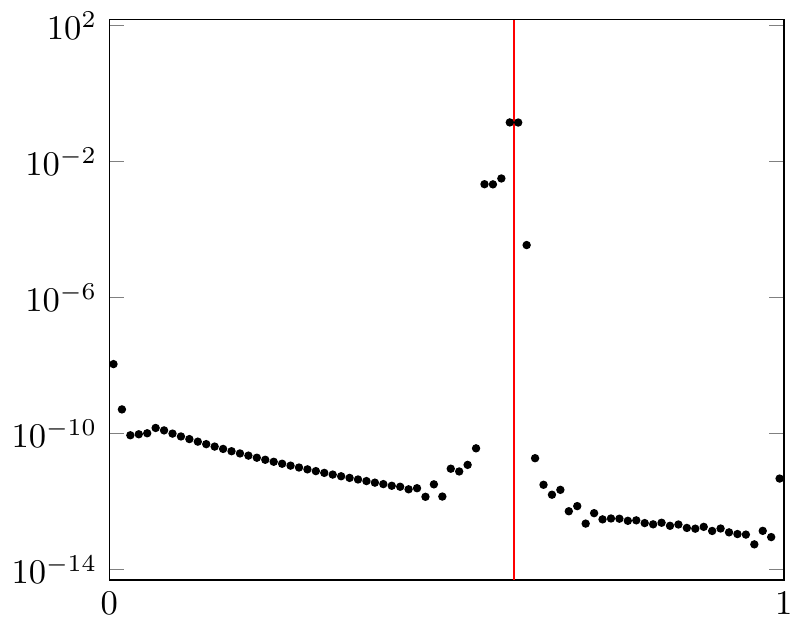}\\\bottomrule
\end{tabular}
\caption{\label{fig:euler_benchmarkI80} Euler example: NL solver for $I=80$ (shock location at an interface).}
\end{figure}

We compare in Fig.~\ref{fig:euler_benchmarkI40} the Static Stencil (first row) used in the original MOOD method versus the Adaptative Stencil (second row) strategies where $I=40$ cells are employed. The simulations, carried out with the non-linear solver, are solely presented since 
essentially equivalent results are obtained with the Time Marching methods. We proceed in the same way with a finer mesh $I=80$ we display in Fig.~\ref{fig:euler_benchmarkI80}. The first column displays the numerical solution for the density variable while the second column gives the CPD together with the stencil configuration by using the color code indicated in this section. No spurious oscillations are observed and the CPD map shows that the Adaptative Stencil method improves the polynomial degrees. We also present in the right column the error cell by cell in order to assess the impact of the stencil shifting close to the discontinuity. We observe a large error cone around the discontinuity where the solution accuracy is significatively reduced due to a spreading of the original error shock. On the contrary, the Adaptative Stencil confines the error in a small region around the shock and we observe a strong improvement of the accuracy in the domains where the solution is smooth.

% $I=40$ (shock location at an interface)

% $I=80$ (shock location at an interface)

To quantify the gain of accuracy, we compute the $L^1$ error over two regions with two meshes $I=40$ and $I=80$. We report in Table~\ref{tab:conv_exE} the errors and convergence rates and show that we recover the optimal sixth-order of accuracy within the smooth area. The Adaptative Stencil turns to be a crucial tool to preserve the solution quality and circumscribe the error propagation deriving from a discontinuity.   
\begin{table}\centering
{\tablepolice
\caption{Errors and convergence order over the regions where the solution of the steady-state Euleur system is smooth. We compute the $L^1$-error on the intervals $[0.050;0.525]$ and $[0.650;0.975]$ while the discontinuity is located at $0.6000$ }
\label{tab:conv_exE}
\begin{tabular}{@{}cccccc@{}}\toprule
& $I$ & $E_{1,\text{L}}[0.050;0.525]$ & ${\mathcal O}_{1,\text{L}}$ & $E_{1,\text{R}}[0.650;0.975]$ & ${\mathcal O}_{1,\text{R}}$\\\midrule
\multirow{2}{*}{MOOD} & 40 & 4.1E$-5$ & --- & 2.3E$-4$ & ---\\
& 80 & 1.9E$-6$ & 4.4 & 1.1E$-5$ & 4.4 \\\midrule
\multirow{2}{*}{AS} & 40 & 1.0E$-9$ & --- & 4.3E$-11$ & ---\\
& 80 & 1.3E$-11$ & 6.3 & 1.2E$-13$ & 8.4 \\\bottomrule
\end{tabular}
}
\end{table}

\begin{figure}\centering
\begin{tabular}{@{}cc@{}c@{}c@{}}\toprule
& exact and numerical solutions & CPD map & errors\\\midrule
MOOD & \includegraphics[width=0.3\textwidth,align=c]{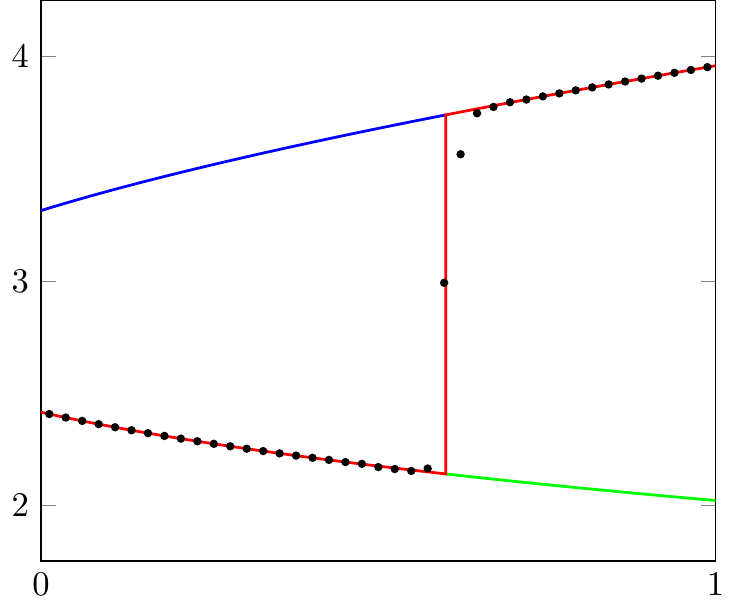} &
\includegraphics[width=0.3\textwidth,align=c]{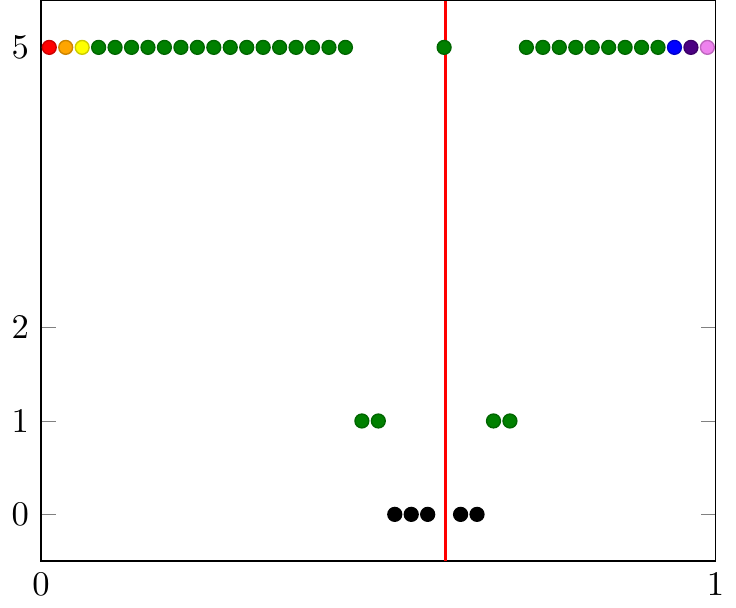} &
\includegraphics[width=0.3\textwidth,align=c]{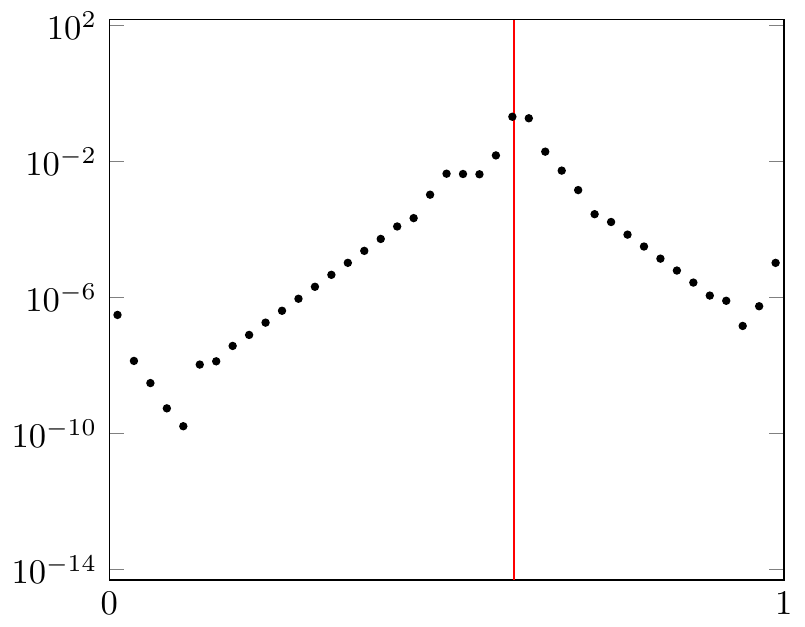}\\
AS & \includegraphics[width=0.3\textwidth,align=c]{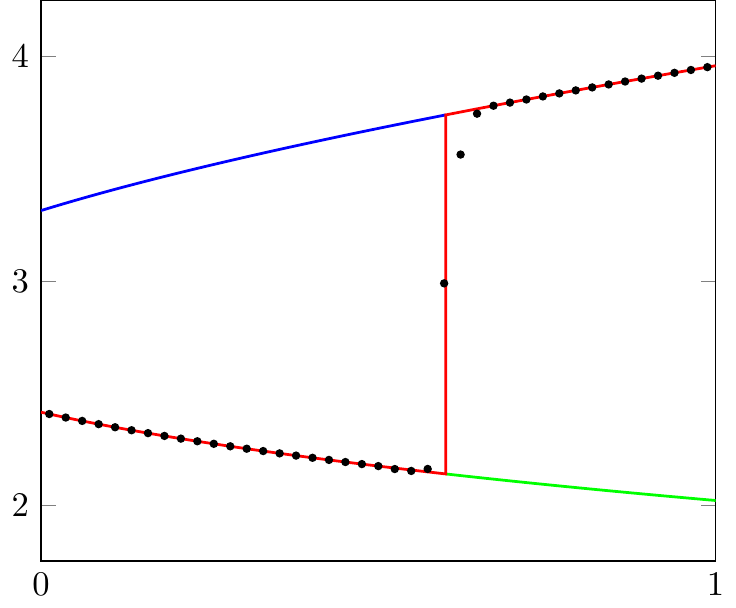} &
\includegraphics[width=0.3\textwidth,align=c]{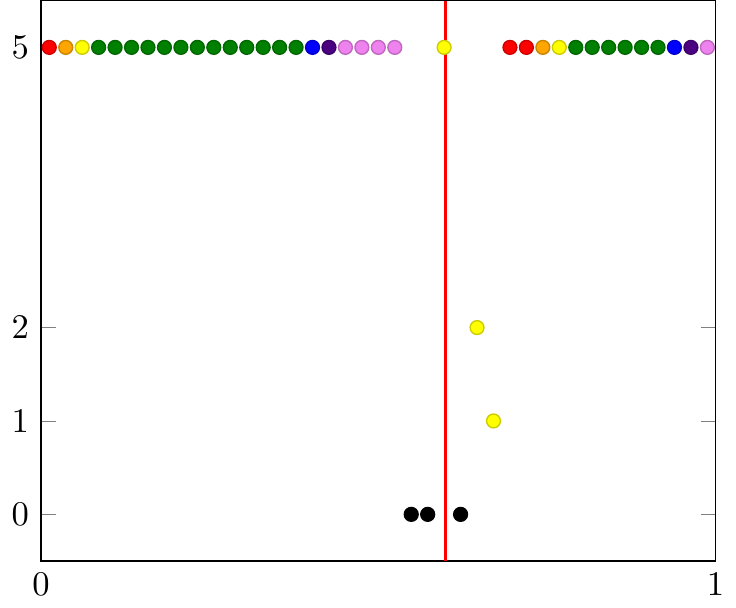} &
\includegraphics[width=0.3\textwidth,align=c]{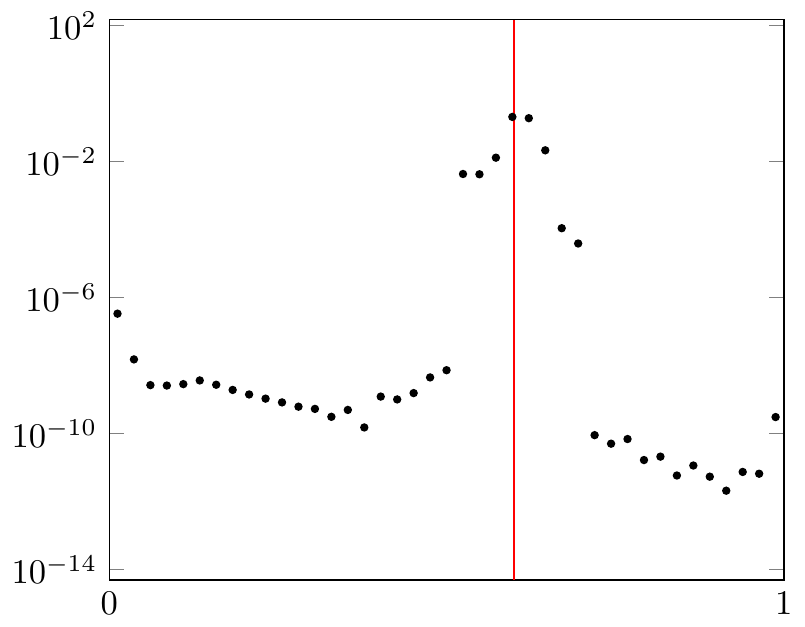}\\\bottomrule
\end{tabular}
\caption{\label{fig:euler_benchmarkI41} Euler example: NL solver for $I=41$ (shock location inside a cell).}
\end{figure}

\begin{figure}\centering
\begin{tabular}{@{}cc@{}c@{}c@{}}\toprule
& exact and numerical solutions & CPD map & errors\\\midrule
MOOD & \includegraphics[width=0.3\textwidth,align=c]{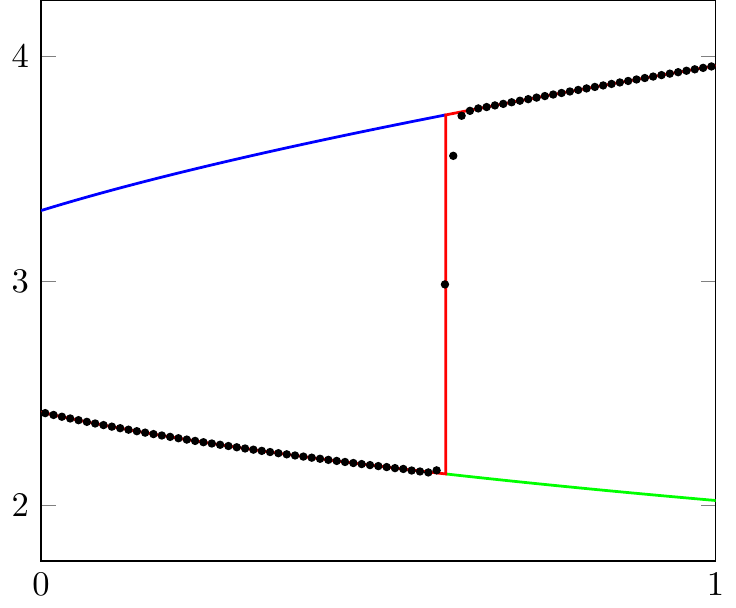} &
\includegraphics[width=0.3\textwidth,align=c]{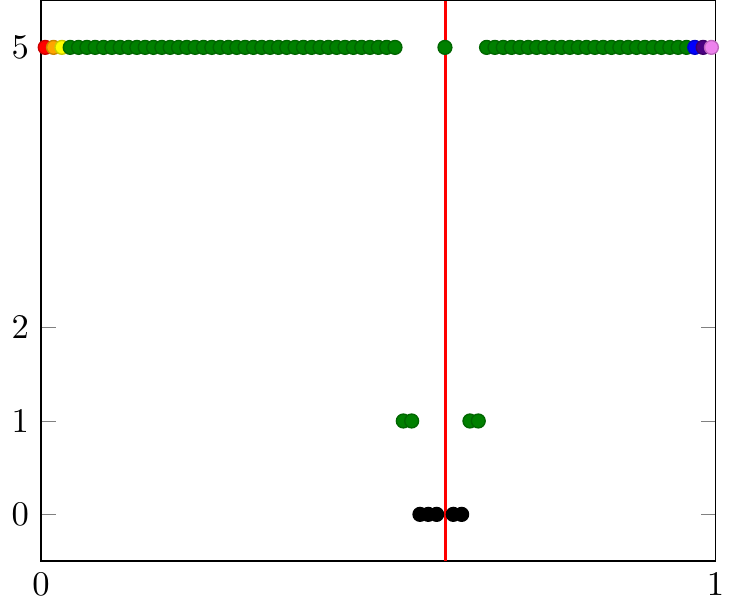} &
\includegraphics[width=0.3\textwidth,align=c]{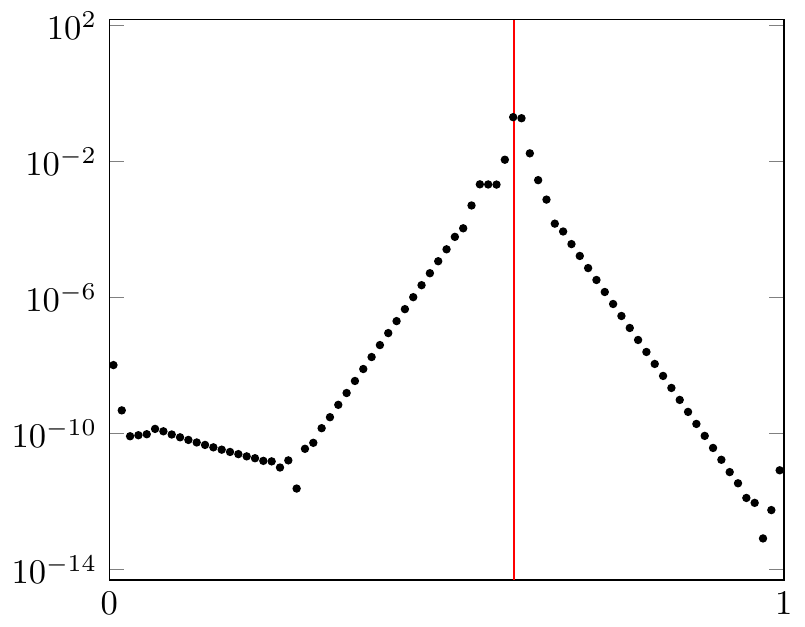}\\
AS & \includegraphics[width=0.3\textwidth,align=c]{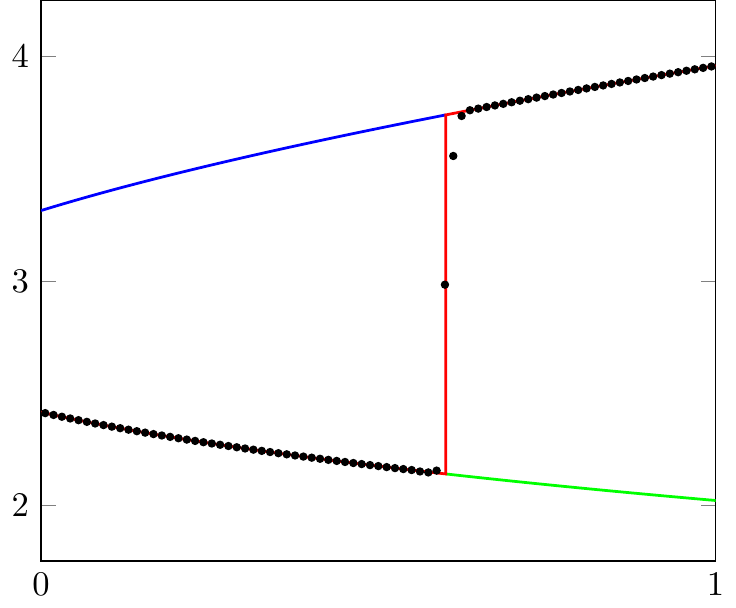} &
\includegraphics[width=0.3\textwidth,align=c]{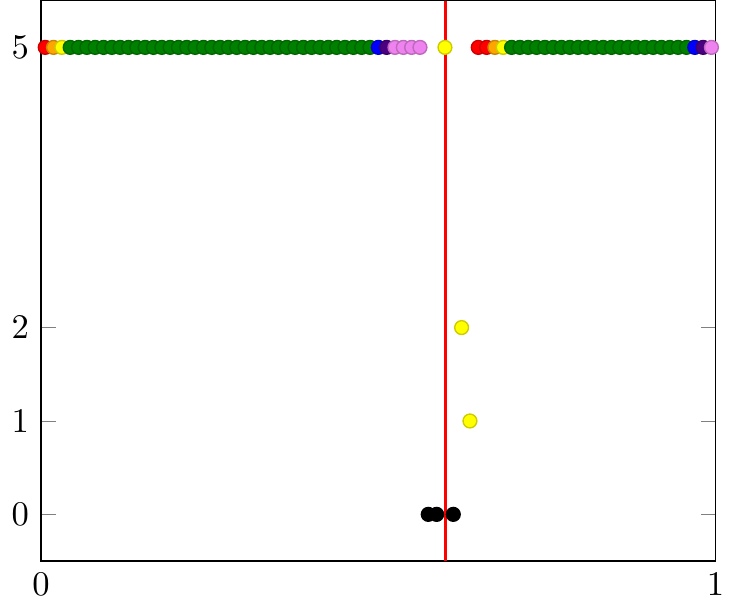} &
\includegraphics[width=0.3\textwidth,align=c]{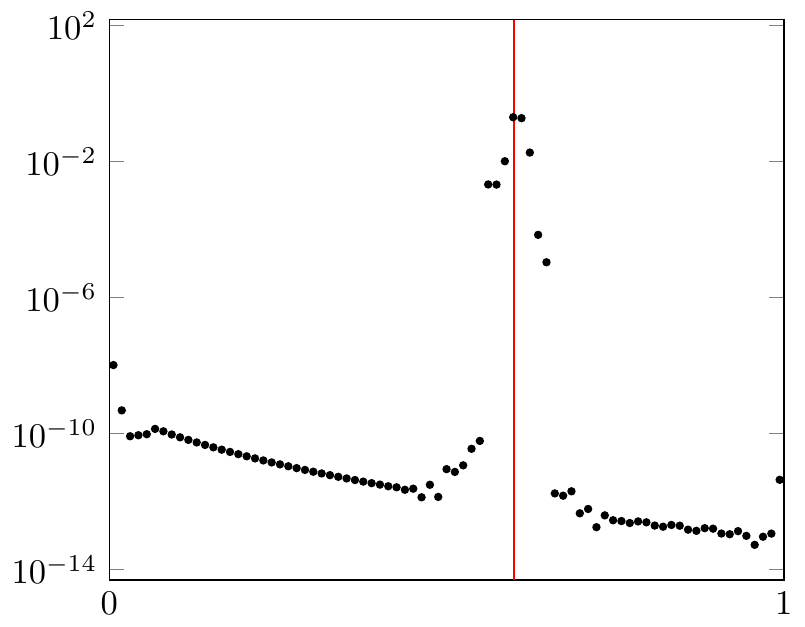}\\\bottomrule
\end{tabular}
\caption{\label{fig:euler_benchmarkI81} Euler example: NL solver for $I=81$ (shock location inside a cell).}
\end{figure}

The shock is located at $x_\text{c}=0.6$ and exactly fit with an inteface between two cells. To check the situation with a shock inside a cell, we repeat the simulation with a slighlty difference by taking $I=41$ and $I=81$ (see Figs.~\ref{fig:euler_benchmarkI41} and \ref{fig:euler_benchmarkI81}). The MOOD method presents a reduction of accuracy of order of one magnitude. On the contrary, we observe that the Adaptative Stencil still performs well and preserve the same accuracy in the domains where the solution is smooth. The method is more robust and less sensitive to the mesh quality. 

% $I=41$ (shock location inside a cell)

% $I=81$ (shock location inside a cell)

% END  NUMERICS
%
%==================================================================================

%==================================================================================
%
%  C O N C L U S I O N
%
\section{Conclusions and perspectives} \label{sec:conclusions}
% Context
In this article we have supplemented the \aposteriori High-Order Finite Volume 
MOOD numerical method with an Adaptive Stencil algorithm. 
Indeed the classical approach only considers
centered stencils for its polynomial reconstruction. However if a discontinuity is present
in the solution, then the reconstruction suffers from the Gibbs phenomenon and some
spurious oscillations may spoil it. 
The classical MOOD approach detects such situations
and reduces the polynomial degree (and \textit{a fortiori} the stencil size) up to
a valid reconstruction is possible or when the degree has dropped to $d=0$, in which 
case the piecewise constant FV data is employed. 
Such a \aposteriori MOOD strategy has revealed fruitful and classical \apriori limiting or stabilization technique can then be replaced by this simple and efficient strategy.
% Why we do this?
However allowing the reconstruction stencil to be adapted, in particular to avoid troubled cells in the vicinity of discontinuity, keeps the polynomial degree at higher value and maintains the accuracy of the overall scheme. 
Remind that for each polynomial degree decrementing, a 
loss in accuracy is generated. In this work we present a proof of concept on 
how an Adaptive Stencil selection can be supplemented to a MOOD scheme for 1D steady-state equations.
First a classical MOOD solution is computed using centered stencils. 
This solution is perfectly acceptable, rather accurate and physically valid. 
Associated to this numerical solution the scheme has appropriately determined 
the Cell Polynomial Degree (CPD) map. Then the stencils are modified and shifted according to the local CPD values, keeping the neighbor cells with high degree and replacing the cells with a low degree. 
Once the new stencils are determined, a MOOD solution is recomputed, starting from the highest degree but using those new reconstruction stencils. 
As such the cost of our MOOD+AS scheme is twice the cost of the classical MOOD scheme.

% Numerical tests
We have shown on 1D advection, B\"urgers', and Euler equations, on smooth and discontinuous solutions
that this AS improves the accuracy of the classical MOOD solution. In particular the accuracy of 
the smooth parts of the flows are often polluted by the presence of a discontinuity, not only in its
vicinity. The AS allows to maintain the highest polynomial degree in many more cells, lowering this
pollution by several orders of magnitude.

% Perspectives
The next step consists in extending this AS to multi-dimensions. In this case the stencil shifting may
occur in several directions that we will have to determine. 
While complex at first glance, this task may be simplified by the large literature on stencil selection for WENO schemes.
Less risky is the extension to time dependent PDEs because we have experimented in this work 
that the time-marching schemes behave appropriately and we are confident
that their adaptation will be relatively easy.

% END  C O N C L U S I O N
%
%==================================================================================

%==================================================================================
% ACKOWLEDGEMENT
\section*{Acknowledgments}
S. Clain and G.J. Machado acknowledge the financial support by FEDER -- Fundo Europeu de Desenvolvimento Regional, through COMPETE 2020 -- Programa Operational Fatores de Competitividade, and the National Funds through FCT -- Funda\c c\~ao para a Ci\^encia e a Tecnologia, project no. UID/FIS/04650/2019.

S. Clain and G.J. Machado acknowledge the financial support by FEDER -- Fundo Europeu de Desenvolvimento Regional, through COMPETE 2020 -- Programa Operacional Fatores de Competitividade, and the National Funds through FCT -- Funda\c c\~ao para a Ci\^encia e a Tecnologia, project no. POCI-01-0145-FEDER-028118.

The material of this research has been partly built and discussed during the SHARK workshops taking place in P\'ovoa de Varzim, Portugal,
\texttt{http://www.SHARK-FV.eu/}.
%==================================================================================

%==================================================================================
% CONFLICT OF INTEREST
\section*{Conflict of interest statement}
On behalf of all authors, the corresponding author states that there is no conflict of interest.

%==================================================================================
% BIBLIOGRAPHY
\bibliographystyle{plain}
\bibliography{biblio}
%==================================================================================

\end{document}